\documentclass[onefignum,onetabnum]{siamart171218}
\usepackage{bbm}
\usepackage{subfig}
\usepackage{amsfonts}
\usepackage{graphicx}
\usepackage{epstopdf}
\usepackage{algorithmic}
\usepackage{ulem}
\DeclareMathOperator*{\argmin}{arg\,min}

\ifpdf
  \DeclareGraphicsExtensions{.eps,.pdf,.png,.jpg}
\else
  \DeclareGraphicsExtensions{.eps}
\fi

\newsiamremark{remark}{Remark}
\newsiamremark{hypothesis}{Hypothesis}
\crefname{hypothesis}{Hypothesis}{Hypotheses}
\newsiamthm{claim}{Claim}
\crefname{claim}{Claim}{Claim}
\newcommand{\alexander}[1]{\textcolor{blue}{\textbf{(alexander)} #1}}

\headers{Variational Extrapolation}{Variational Extrapolation}

\title{Variational Extrapolation of Implicit Schemes\\ for General Gradient Flows\thanks{\textbf{Funding:} Alexander Zaitzeff and Selim Esedo\=glu gratefully acknowledge support from the NSF grant DMS-1719727. Krishnna Garikipati acknowledges NSF grants DMREF-1436154 and DMREF-1729166.}}

\author{Alexander Zaitzeff\thanks{Department of Mathematics, University of Michigan, Ann Arbor, MI 48109, USA
  (\email{azaitzef@umich.edu}, \email{esedoglu@umich.edu}).}
\and Selim Esedo\=glu\footnotemark[2]
\and Krishna Garikipati\thanks{Departments of Mechanical Engineering, and Mathematics, Michigan Institute for Computational Discovery \& Engineering, University of Michigan, Ann Arbor, MI 48109, USA (\email{krishna@umich.edu}).} }

\usepackage{amsopn}

\ifpdf
\hypersetup{
  pdftitle={Variational Extrapolation},
  pdfauthor={Zaitzeff, Esedoglu, Garikipati}
}
\fi

\begin{document}

\maketitle


\begin{abstract}
We introduce a class of unconditionally energy stable, high order accurate schemes for gradient flows in a very general setting.
The new schemes are a high order analogue of the minimizing movements approach for generating a time discrete approximation to a gradient flow by solving a sequence of optimization problems.
In particular, each step entails minimizing the associated energy of the gradient flow plus a movement limiter term that is, in the classical context of steepest descent with respect to an inner product, simply quadratic.
A variety of existing unconditionally stable numerical methods can be recognized as (typically just first order accurate in time) minimizing movement schemes for their associated evolution equations, already requiring the optimization of the energy plus a quadratic term at every time step. Therefore, our approach gives a painless way to extend these to high order accurate in time schemes while maintaining their unconditional stability. In this sense, it can be viewed as a variational analogue of Richardson extrapolation.    
\end{abstract}

\begin{keywords}
Extrapolation, Gradient flows, High order schemes, Unconditional stability, Minimizing movements
\end{keywords}

\begin{AMS}
   65M12, 65K10, 65L06, 65L20
\end{AMS}

\section{Introduction}
We are concerned with numerical schemes for evolution equations that arise as gradient flow (steepest descent) for an energy $E:H \rightarrow \mathbb{R}$, where $H$ is a Hilbert space with inner product $\langle\cdot,\cdot\rangle$:
\begin{equation}
    \label{eq:de}
    u'=-\nabla_H E(u).
\end{equation}
Equation \cref{eq:de} may represent a (scalar or vectorial) ordinary or partial differential equation. 
A fundamental property of equation \cref{eq:de} is that it dissipates the energy over time:
\begin{equation*}
\frac{d}{dt} E(u)=\langle \nabla_H E(u),u'\rangle=-||\nabla_H E(u)||^2\leq 0.
\end{equation*}
Our focus is on unconditionally energy stable, high order in time discretizations. To be precise, by energy stable we mean the following dissipative property: 
\begin{equation}
\label{eq:energymin}
    E(u_{n+1})\leq E(u_n)
\end{equation}
where $u_n$ denotes the approximation to the solution at the $n$-th time step.
Thus, in the context of PDEs, where $H$ is infinite dimensional, we are concerned with discrete in time, continuous in space schemes.

The backward Euler method for the abstract equation (\ref{eq:de}), with time step size $k>0$, reads
\begin{equation}
\label{eq:bE}
\frac{u_{n+1} - u_n}{k} = -\nabla_H E(u_{n+1}).
\end{equation}
As is well known and immediate to see, a solution $u_{n+1}$ for the implicit scheme (\ref{eq:bE}) can be found via the following optimization problem
\begin{equation}
\label{eq:minmov}
u_{n+1}=\argmin_u \bigg(E(u)+\frac{1}{2k}||u-u_n||^2\bigg)
\end{equation}
since (\ref{eq:bE}) is the Euler-Lagrange equation for the optimization (\ref{eq:minmov}); here, $\| \cdot \|^2 = \langle \cdot , \cdot \rangle$.
It follows that 
\begin{equation}
E(u_{n+1})\leq E(u_{n+1})+\frac{1}{2k}||u_{n+1}-u_n||^2\leq E(u_{n})+\frac{1}{2k}||u_n-u_n||^2=E(u_{n})
\end{equation}
so that scheme (\ref{eq:bE}) is unconditionally stable, provided that optimization problem (\ref{eq:minmov}) can be solved.

Energetic formulation (\ref{eq:minmov}) of the backward Euler scheme (\ref{eq:bE}) is often referred to as {\it minimizing movements}.
It enables extending numerical schemes for the stationary optimization problem $\min_u E(u)$ to the dynamic, evolutionary problem (\ref{eq:de}) provided an additional, typically quadratic term in the cost function can be accommodated.
The quadratic term $\frac{1}{2k} \| u-u_n \|^2$ in (\ref{eq:minmov}) is often referred to as the {\it movement limiter}, as it opposes deviation from the current configuration $u_n$.
It encodes the inner product with respect to which the gradient flow is being generated.
Beyond numerical analysis and computation, minimizing movements approximation of gradient flows have been instrumental in the analysis of evolution equations of the form (\ref{eq:de}), e.g. in defining and finding weak solutions beyond the formation of singularities when classical notions of solution cease to exist.

\smallskip

The following combination of desirable properties distinguish the new schemes introduced in this paper:
\begin{enumerate}
\item Complete generality. There is no assumption (e.g. convexity) on the energy $E$ in (\ref{eq:de}) beyond sufficient differentiability.
\item Unconditional energy stability.
\item High (at least up to third) order accuracy.
\item Each time step requires a few standard minimizing movements solves, equivalent to backward Euler substeps, or optimization of the associated energy plus a quadratic term.
\end{enumerate}
\smallskip

\noindent Property 4 is perhaps the most unique and appealing aspect of the new framework: There are many existing schemes that can be recognized as some form of minimizing movements, sometimes relying on efficient optimization algorithms to solve (\ref{eq:bE}) via (\ref{eq:minmov}).
Our contribution shows how to painlessly jack up the order of accuracy of these schemes while preserving unconditional stability, relying only on a black-box implementation of the standard backward Euler scheme.
In that sense, our new schemes can be understood as a variational analogue of Richardson extrapolation on \cref{eq:bE}, which in its standard form lacks the stability guarantees of our new schemes.

Many general purpose numerical schemes can certainly be used for solving \cref{eq:de}, such as multistep or Runge-Kutta methods \cite{butcher2016numerical}.
However, the energy stability of the standard examples of such schemes is either not immediate, or not true at all, at the level of generality we seek here, when applied to an equation of the form \cref{eq:de}.
Our focus is on high order schemes whose stability can be guaranteed over an entire class of evolution laws, namely gradient flows \cref{eq:de}.
Nevertheless, after some appropriate transformations, the new schemes we propose can be seen as a new, special class of diagonally implicit Runge-Kutta (DIRK) schemes tailored to these important dynamics.
In the extensive literature on Runge-Kutta methods, one of the related contributions to the nonlinear notion of stability \cref{eq:energymin} we seek is B-stability for evolutions that satisfy a monotonicity (contractivity) condition \cite{butcher1975}.
In the context of gradient flows, this requires convexity of the energy $E$ in \cref{eq:de}, which is too restrictive for the applications we have in mind (see e.g. Examples \cref{eq:ac} and \cref{eq:ch} in Section \ref{subsec:PDE}).
Very recently, \cite{shin2017unconditionally} \& \cite{shin2020energy} propose high order Runge-Kutta schemes for gradient flows with stability guarantees. Among these, \cite{shin2020energy} concerns fully implicit schemes, as in the present work, but is again restricted to convex energies as in earlier works on B-convexity. The paper \cite{shin2017unconditionally} studies implicit-explicit schemes that, in the spirit of convexity splitting \cite{eyre1998unconditionally}, break up the energy into convex and concave parts, and treat the convex part implicitly and the concave part explicitly.
The present work differs in placing no convexity assumptions on $E$, which is treated fully implicitly.
An example where the energy is in fact {\it concave}, yet the optimization \cref{eq:minmov} is solvable at very low cost, is the threshold dynamics algorithm for motion by mean curvature \cite{mbo92,mbo94} that is known to be unconditionally energy stable \cite{eo}. We show in \cite{alex2019second} how ideas developed in the present paper can be used to jack up the order of accuracy of this intriguing algorithm while preserving its desirable stability properties, which appears to be beyond the scope of previous contributions. 
See also Remark \ref{remark:thresh} of Section \ref{subsec:PDE} in this context.
Finally, we also mention recent work on the scalar auxiliary variable method ~\cite{shen2018scalar} as another approach focusing on unconditional energy stability for gradient flows.

The rest of the paper is organized as follows:
\begin{itemize}
\item \Cref{sec:gf} presents the general framework for the new scheme, focusing on unconditional energy stability.
\item \Cref{sec:ex} focuses on consistency, showing how to attain 2nd and 3rd order accuracy.
 \item \Cref{sec:examples} gives 2nd and 3rd order examples of the new schemes.
\item \Cref{sec:nr} presents numerical convergence studies on a number of well-known ordinary and partial differential equations that are gradient flows.
\end{itemize}
The code for  \cref{sec:nr} is publicly available, and can be found at \url{https://github.com/AZaitzeff/gradientflow}.

\section{The New Schemes: Stability}
\label{sec:gf}

In this section, we formulate a wide class of numerical schemes that are energy stable by construction.
We thus place stability front and center, leaving consistency to be dealt with subsequently.
It is therefore important to allow many degrees of freedom in the scheme at this stage, in the form of a large number of coefficients, that will eventually be chosen, in the next section, to attain consistency at a high order of accuracy.

Our method is a linear $M$-stage scheme of the following form:
\begin{enumerate}
    \item Set $U_0 = u_n$
    \item For $m=1,\ldots,M$:
\begin{equation}
\label{eq:ms}
U_m=\argmin_u \bigg(E(u)+\sum^{m-1}_{i=0}\frac{\gamma_{m,i}}{2k}||u-U_i||^2\bigg).
\end{equation}
\item Set $u_{n+1}=U_M$
\end{enumerate}

Notice that the proposed scheme \cref{eq:ms}, as promised, merely requires the solution of exactly the same type of problem at every time step as the standard backward Euler scheme: minimization of the associated energy plus a quadratic term.

At this point, it is not clear why a scheme such as \cref{eq:ms} should dissipate energy $E$ at every iteration as in \cref{eq:energymin}. However, in this section we establish quite broad conditions on the coefficients $\gamma_{m,i}$ that ensure energy dissipation \cref{eq:energymin}; this is the essential observation at the heart of the present paper. 
To demonstrate the idea, consider the following two-stage special case of scheme \cref{eq:ms}:
\begin{align}
\label{eq:twostage1}
&U_1=\argmin_u \bigg(E(u) +\frac{\gamma_{1,0}}{2k}||u-u_n||^2\bigg)\\
\label{eq:twostage2}
&u_{n+1}=\argmin_u \bigg(E(u) +\frac{\gamma_{2,0}}{2k}||u-u_n||^2+\frac{\gamma_{2,1}}{2k}||u-U_1||^2\bigg)
\end{align}
and impose the conditions
\begin{equation}
\label{eq:twostage3} 
\gamma_{1,0}-\frac{\gamma_{2,0}^2}{\gamma_{2,0}+\gamma_{2,1}} \geq 0 \mbox{ and } \gamma_{2,0}+\gamma_{2,1}>0
\end{equation}
on the parameters.
Set $\theta = \frac{\gamma_{2,0}}{\gamma_{2,0}+\gamma_{2,1}}$. 
Note that \eqref{eq:twostage2} is equivalent to
\begin{equation}
\label{eq:twostage4}
u_{n+1}  = \argmin_u \bigg( E(u)+\frac{\gamma_{2,0}+\gamma_{2,1}}{2k} \left\| u - \big( \theta u_n + (1-\theta) U_1 \big) \right\|^2 \bigg).
\end{equation}
This can be seen by expanding the norm squared and comparing
the quadratic and linear terms in $u$. The constant terms are not equal but
that does not matter for the minimization.\\
We have
\small
\begin{align*}
E(u_{n+1}) &\leq E(u_{n+1}) + \frac{\gamma_{2,0}+\gamma_{2,1}}{2k} \left\| u_{n+1} - \big( \theta u_n + (1-\theta) U_1 \big) \right\|^2 && \text{(by \eqref{eq:twostage3})}\\
&\leq E(U_1) + \frac{\gamma_{2,0}+\gamma_{2,1}}{2k} \left\|U_1 - \big( \theta u_n + (1-\theta) U_1 \big) \right\|^2 && \text{(by \eqref{eq:twostage4})}\\
&= E(U_1) + \frac{\gamma_{2,0}^2}{(\gamma_{2,1}+\gamma_{2,0})2k} \left\| U_1 - u_n \right\|^2\\
&\leq   E(U_1) + \frac{\gamma_{1,0}}{2k} \left\| U_1 - u_n \right\|^2 && \text{(by \eqref{eq:twostage3})}\\
&\leq E(u_n) && \text{(by \eqref{eq:twostage1})}
\end{align*}
\normalsize
establishing unconditional energy stability of scheme (\ref{eq:twostage1}) \& (\ref{eq:twostage2}) under the condition (\ref{eq:twostage3}) on its parameters. We offer some insight to the conditions in  \cref{eq:twostage3}. First, the condition $\gamma_{2,0}+\gamma_{2,1}>0$ is reasonable as it requires that the function being minimized in \cref{eq:twostage2} goes to $+\infty$ as $||u||\to \infty$. What is more surprising is that one of the second stage coefficients can be negative while maintaining unconditionally stability. We can `reward' the distance to one of the previous stages as long as the distance to the other stage is penalized sufficiently strongly. The condition $\gamma_{1,0} \geq \frac{\gamma_{2,0}^2}{\gamma_{2,0}+\gamma_{2,1}}$ requires that the penalization in the first stage has to be strong relative to the penalization in the second stage.\\
We will now extend this discussion to the general, $M$-stage case of scheme \cref{eq:ms}:

\begin{theorem}
\label{claim:ms}
Define the following auxiliary quantities in terms of the coefficients $\gamma_{m,i}$ of scheme \cref{eq:ms}:
\begin{align}
\label{eq:tildegamma}
&\tilde{\gamma}_{m,i}=\gamma_{m,i}-\sum_{j=m+1}^M\tilde{\gamma}_{j,i}\frac{\tilde{S}_{j,m}}{\tilde{S}_{j,j}}\\
\label{eq:S}
&\tilde{S}_{j,m}=\sum_{i=0}^{m-1} \tilde{\gamma}_{j,i}
\end{align}
Where if $m=M$, the sum in \cref{eq:tildegamma} is considered to be zero. If $\tilde{S}_{m,m}>0$ for $m=1,\ldots,M$, then scheme \cref{eq:ms} satisfies the energy stability condition (\ref{eq:energymin}): For every $n=0,1,2,\ldots$ we have $E(u_{n+1}) \leq E(u_n)$.
\end{theorem}
As we will see in Section \ref{sec:ex}, the conditions on the parameters $\gamma_{i,j}$ of scheme \cref{eq:ms} imposed in \cref{claim:ms} are loose enough to enable meeting consistency conditions to high order.
We will establish \cref{claim:ms} with the help of the following two lemmas.The first lemma is the multi-step version of the equivalence of \cref{eq:twostage2} and \cref{eq:twostage4} in our two step example:

\begin{lemma}
Let the auxiliary quantities $\tilde{S}_{j,m}$, and $\tilde{\gamma}_{m,i}$ be defined as in \cref{claim:ms}. 
We have
\label{lem:equiv}
\begin{align*}
&\argmin_u \bigg( E(u)+\sum^{m-1}_{i=0}\frac{\gamma_{m,i}}{2k}||u-U_i||^2\bigg)\\=
&\argmin_u \bigg(E(u)+\frac{1}{2k}\sum_{j=m}^M \frac{\tilde{S}_{j,m}^2}{\tilde{S}_{j,j}} ||u-\sum_{i=0}^{m-1} \frac{\tilde{\gamma}_{j,i}}{\tilde{S}_{j,m}} U_i||^2\bigg)
\end{align*}
\end{lemma}

\begin{proof}
As in the two step case the proof consists of expanding the norm squared terms and showing that all the quadratic and linear terms of $u$ are equal. First the expansion of $\sum^{m-1}_{i=0}\frac{\gamma_{m,i}}{2k}||u-U_i||^2$ is 
\begin{align}
\label{eq:firstexp}
\frac{||u||^2}{2k}\sum^{m-1}_{i=0}\gamma_{m,i} - \frac{1}{k}\langle u,\sum^{m-1}_{i=0}\gamma_{m,i}U_i\rangle+\text{terms that do not depend on $u$.}
\end{align}
Next, we will establish two identities to help us expand
\small\[\frac{1}{2k}\sum_{j=m}^M \frac{\tilde{S}_{j,m}^2}{\tilde{S}_{j,j}} ||u-\sum_{i=0}^{m-1} \frac{\tilde{\gamma}_{j,i}}{\tilde{S}_{j,m}} U_i||^2.\]
\normalsize
First by rearranging \cref{eq:tildegamma},
\small
\begin{equation}
\label{eq:gammaind}
 \gamma_{m,i}=\sum_{j=m}^M\tilde{\gamma}_{j,i}\frac{\tilde{S}_{j,m}}{\tilde{S}_{j,j}}.
 \end{equation}
 \normalsize
Next, an identity of $\tilde{S}_{m,m}$:
\small
\begin{align*}
\tilde{S}_{m,m}&=\sum_{i=0}^{m-1} \tilde{\gamma}_{m,i}=\sum_{i=0}^{m-1}\bigg[ \gamma_{m,i}-\sum_{j=m+1}^M\tilde{\gamma}_{j,i}\frac{\tilde{S}_{j,m}}{\tilde{S}_{j,j}}\bigg]\\&= \sum_{i=0}^{m-1} \gamma_{m,i}-\sum_{j=m+1}^M\bigg[\sum_{i=0}^{m-1}\tilde{\gamma}_{j,i}\bigg]\frac{\tilde{S}_{j,m}}{\tilde{S}_{j,j}}=\sum_{i=0}^{m-1} \gamma_{m,i}-\sum_{j=m+1}^M \frac{\tilde{S}_{j,m}^2}{\tilde{S}_{j,j}}.
\end{align*}
\normalsize
We use this identity to establish the following:
\small
\begin{equation}
\label{eq:idenitysum}
    \sum_{j=m}^M \frac{\tilde{S}_{j,m}^2}{\tilde{S}_{j,j}}=\tilde{S}_{m,m}+\sum_{j=m+1}^M \frac{\tilde{S}_{j,m}^2}{\tilde{S}_{j,j}}=\sum_{i=0}^{m-1} \gamma_{m,i}-\sum_{j=m+1}^M \frac{\tilde{S}_{j,m}^2}{\tilde{S}_{j,j}}+\sum_{j=m+1}^M \frac{\tilde{S}_{j,m}^2}{\tilde{S}_{j,j}}=\sum_{i=0}^{m-1} \gamma_{m,i}.
\end{equation}

\normalsize
Now we can calculate the expansion:
\small
\begin{align*}
&\frac{1}{2k}\sum_{j=m}^M \frac{\tilde{S}_{j,m}^2}{\tilde{S}_{j,j}} ||u-\sum_{i=0}^{m-1} \frac{\tilde{\gamma}_{j,i}}{\tilde{S}_{j,m}} U_i||^2\\
=&\frac{||u||^2}{2k} \sum_{j=m}^M \frac{\tilde{S}_{j,m}^2}{\tilde{S}_{j,j}} -\frac{1}{k}\langle u,\sum_{i=0}^{m-1} \sum_{j=m}^M \tilde{\gamma}_{j,i} \frac{\tilde{S}_{j,m}}{\tilde{S}_{j,j}} U_i \rangle+\text{terms that do not depend on $u$}\\
=&\frac{||u||^2}{2k}\sum^{m-1}_{i=0}\gamma_{m,i} - \frac{1}{k}\langle u,\sum^{m-1}_{i=0}\gamma_{m,i}U_i\rangle+\text{terms that do not depend on $u$.}
\end{align*}\\
\normalsize Where the last equality follows from \cref{eq:gammaind} and \cref{eq:idenitysum}. Since this expansion matches \cref{eq:firstexp} up to a constant in $u$ the proof is complete.
\end{proof}

Now we will use \cref{lem:equiv} to relate the energy of sub-step $m$  to sub-step $m-1$. This lemma is the crux of the proof of the theorem and where we use the condition that $\tilde{S}_{m,m}>0$ for all $m$.

\begin{lemma}
\label{lem:step}
Let the auxiliary quantities $\tilde{S}_{j,m}$, $\tilde{\gamma}_{m,i}$ be given in \cref{claim:ms} and let $\tilde{S}_{m,m}>0$ for $m=1,\ldots,M$. Then
\begin{align*}
    &E(U_m)+\frac{1}{2k}\sum_{j=m}^M \frac{\tilde{S}_{j,m}^2}{\tilde{S}_{j,j}} ||U_m-\sum_{i=0}^{m-1} \frac{\tilde{\gamma}_{j,i}}{\tilde{S}_{j,m}} U_i||^2\\\leq & E(U_{m-1})+\frac{1}{2k}\sum_{j=m-1}^M \frac{\tilde{S}_{j,m-1}^2}{\tilde{S}_{j,j}} ||U_{m-1}-\sum_{i=0}^{m-2} \frac{\tilde{\gamma}_{j,i}}{\tilde{S}_{j,m-1}} U_i||^2
\end{align*}
\end{lemma}

\begin{proof}
\small
By \cref{eq:ms} \& \cref{lem:equiv},
\[
U_m=\argmin_u E(u)+\frac{1}{2k}\sum_{j=m}^M \frac{\tilde{S}_{j,m}^2}{\tilde{S}_{j,j}} ||u-\sum_{i=0}^{m-1} \frac{\tilde{\gamma}_{j,i}}{\tilde{S}_{j,m}} U_i||^2.
\]
\normalsize
Since $U_m$ is the minimizer of the above optimization problem
\small
\begin{align}
\label{eq:lem2step0}
&E(U_m)+\frac{1}{2k}\sum_{j=m}^M \frac{\tilde{S}_{j,m}^2}{\tilde{S}_{j,j}} ||U_m-\sum_{i=0}^{m-1} \frac{\tilde{\gamma}_{j,i}}{\tilde{S}_{j,m}} U_i||^2\\
\label{eq:lem2step1}
&\leq E(U_{m-1})+\frac{1}{2k}\sum_{j=m}^M \frac{\tilde{S}_{j,m}^2}{\tilde{S}_{j,j}} ||U_{m-1}-\sum_{i=0}^{m-1} \frac{\tilde{\gamma}_{j,i}}{\tilde{S}_{j,m}} U_i||^2
\end{align}
\normalsize
Next using the definition of auxiliary variables we can state an identity that will simplify \cref{eq:lem2step1}. For $m>1$ and $j\geq m$ 
\small
\begin{align*}
&\frac{\tilde{S}_{j,m}^2}{\tilde{S}_{j,j}} ||U_{m-1}-\sum_{i=0}^{m-1} \frac{\tilde{\gamma}_{j,i}}{\tilde{S}_{j,m}} U_i||^2=\frac{\tilde{S}_{j,m}^2}{\tilde{S}_{j,j}} ||U_{m-1}\bigg(1-\frac{\tilde{\gamma}_{j,m-1}}{\tilde{S}_{j,m}}\bigg)-\sum_{i=0}^{m-2} \frac{\tilde{\gamma}_{j,i}}{\tilde{S}_{j,m}} U_i||^2\\
=&\frac{\tilde{S}_{j,m}^2}{\tilde{S}_{j,j}} ||U_{m-1}\bigg(\frac{\tilde{S}_{j,m-1}}{\tilde{S}_{j,m}}\bigg)-\sum_{i=0}^{m-2} \frac{\tilde{\gamma}_{j,i}}{\tilde{S}_{j,m}} U_i||^2=\frac{\tilde{S}_{j,m-1}^2}{\tilde{S}_{j,j}} ||U_{m-1}-\sum_{i=0}^{m-2} \frac{\tilde{\gamma}_{j,i}}{\tilde{S}_{j,m-1}} U_i||^2.
\end{align*}
\normalsize
Using this identity \cref{eq:lem2step1} is equal to
\small
\begin{equation}
\label{eq:lem2step2}
E(U_{m-1})+\frac{1}{2k}\sum_{j=m}^M \frac{\tilde{S}_{j,m-1}^2}{\tilde{S}_{j,j}} ||U_{m-1}-\sum_{i=0}^{m-2} \frac{\tilde{\gamma}_{j,i}}{\tilde{S}_{j,m-1}} U_i||^2
\end{equation}
\normalsize
Now since $\tilde{S}_{m-1,m-1}>0$,
\small
\begin{equation}
   \label{eq:lem2step3} 
    \frac{\tilde{S}_{m-1,m-1}^2}{\tilde{S}_{m-1,m-1}} ||U_{m-1}-\sum_{i=0}^{m-2} \frac{\tilde{\gamma}_{m-1,i}}{\tilde{S}_{m-1,m-1}} U_i||^2>0.
\end{equation}
\normalsize
By adding \cref{eq:lem2step3} to \cref{eq:lem2step2}, we have that \cref{eq:lem2step2} is less than or equal to
\small
\begin{align*}
E(U_{m-1})+\frac{1}{2k}\sum_{j=m-1}^M \frac{\tilde{S}_{j,m-1}^2}{\tilde{S}_{j,j}} ||U_{m-1}-\sum_{i=0}^{m-2} \frac{\tilde{\gamma}_{j,i}}{\tilde{S}_{j,m-1}} U_i||^2
\end{align*}
\normalsize
concluding the proof.
\end{proof}

\begin{proof}(of theorem)
The main idea of the proof is to use \cref{lem:step} repeatedly to relate the energy of $E(u_{n+1})$ to $E(u_n)$. First, since $\tilde{S}_{M,M}>0$
\[E(u_{n+1})=E(U_M)
\leq E(U_M)+\frac{1}{2k} \frac{\tilde{S}_{M,M}^2}{\tilde{S}_{M,M}} ||U_M-\sum_{i=0}^{M-1} \frac{\tilde{\gamma}_{M,i}}{\tilde{S}_{M,M}} U_i||^2\]
The right hand side of the equation is of the form required by \cref{lem:step}. By using the \cref{lem:step} repeatedly we have
\small
\begin{align*}
&E(U_M)+\frac{1}{2k} \frac{\tilde{S}_{M,M}^2}{\tilde{S}_{M,M}} ||U_M-\sum_{i=0}^{M-1} \frac{\tilde{\gamma}_{M,i}}{\tilde{S}_{M,M}} U_i||^2\\
\leq& E(U_{M-1})+\frac{1}{2k}\sum_{j=M-1}^M \frac{\tilde{S}_{j,M-1}^2}{\tilde{S}_{j,j}} ||U_{M-1}-\sum_{i=0}^{M-2} \frac{\tilde{\gamma}_{j,i}}{\tilde{S}_{j,M-1}} U_i||^2\\
& \vdots \\\leq& E(U_1)+\frac{1}{2k}\sum_{j=1}^M \frac{\tilde{S}_{j,1}^2}{\tilde{S}_{j,j}} ||U_1- \frac{\tilde{\gamma}_{j,0}}{\tilde{S}_{j,1}} U_0||^2.
\end{align*}
\normalsize 
By \cref{eq:ms} and \cref{lem:equiv}
\small
\[U_1=\argmin_u E(u)+\frac{1}{2k}\sum_{j=1}^M \frac{\tilde{S}_{j,1}^2}{\tilde{S}_{j,j}} ||u- \frac{\tilde{\gamma}_{j,0}}{\tilde{S}_{j,1}} U_0||^2\]
\normalsize
so
\small
\[E(U_1)+\frac{1}{2k}\sum_{j=1}^M \frac{\tilde{S}_{j,1}^2}{\tilde{S}_{j,j}} ||U_1- \frac{\tilde{\gamma}_{j,0}}{\tilde{S}_{j,1}} U_0||^2 \leq E(U_0)+\frac{1}{2k}\sum_{j=1}^M \frac{\tilde{S}_{j,1}^2}{\tilde{S}_{j,j}} ||U_0-U_0||^2=E(u_n)\]
\normalsize
completing the proof of the theorem.
\end{proof}

Now the condition that $\tilde{S}_{m,m}>0$ for $m=1,\ldots,M$ is the multi-step equivalent of \cref{eq:twostage3} in the two step case. Given the $\gamma$'s you can calculate the auxiliary quantities \cref{eq:tildegamma} and \cref{eq:S}  explicitly as follows:\\
for $m=M,M-1,\ldots,1$:
\begin{enumerate}
    \item Calculate $\tilde{\gamma}_{m,i}$ for $i=1,2,\ldots,m$
    \item Calculate $\tilde{S}_{j,m}$ for $j=m,m+1,\ldots,M$.
\end{enumerate}
Thus given $\gamma$'s we can easily check if they satisfy the hypothesis of \cref{claim:ms}.

\section{The New Schemes: Consistency}
\label{sec:ex}
We now turn to the question of whether the coefficients $\gamma_{m,i}$ in scheme \cref{eq:ms} can be chosen to ensure its high order consistency with the abstract evolution law \cref{eq:de}. As mentioned before, the schemes are diagonal implicit Runge-Kutta, whose order conditions are well established (for example in ~\cite{butcher2016numerical}). For completeness, we derive the conditions here.
From \cref{eq:ms}, each stage $U_m$ satisfies the Euler-Lagrange equation:
\begin{equation}
    \label{eq:ms2}
    \bigg[\sum^{m-1}_{i=0}\gamma_{m,i}\bigg]U_{m}+k\nabla_HE(U_m)=\sum^{m-1}_{i=0}\gamma_{m,i}U_i.
\end{equation}
The consistency equations for the $\gamma$s are found by  carrying out a Taylor series expansion of $U_m$ around $U_0=u(t_0)$. We will calculate the one-step error. For $n\in\{1,2,3,\ldots\}$, let $D^n E(u) : H^n \to \mathbb{R}$ denote the multilinear form given by
$$ D^n E(u)(v_1,\ldots,v_n) = \left. \frac{\partial^n}{\partial s_1 \cdots \partial s_n} E(u+s_1 v_1 + s_2 v_2 + \cdots + s_n v_n) \right|_{s_1 = s_2 = \cdots = s_n = 0}$$
so that the linear functional $D^n E(u)(v_1,v_2,\ldots,v_{n-1},\cdot) : H \to \mathbb{R}$ may be identified with an element of $H$, which will be denoted simply as $D^n E(v_1,v_2,\ldots,v_{n-1})$ in what follows.
We begin with the exact solution starting from $u(t_0)$:
\[
\begin{cases}
u_t=-\nabla E(u) & t>t_0 \\
u(t_0)=U_0
  \end{cases}
\]
The Taylor expansion of $u(k+t_0)$ around $t_0$ is
\small
\begin{align}
\begin{split}
\label{eq:true}
u(k+t_0)=&u(t_0)+ku_t(t_0)+\frac{1}{2}k^2u_{tt}(t_0)+\frac{1}{6}k^3u_{ttt}(t_0)+\text{h.o.t.}\\=&U_0-k DE(U_0)+\frac{1}{2}k^2D^2E(U_0)DE(U_0)\\&-\frac{1}{6}k^3\big[D^2E(U_0)\left(D^2E(U_0)\left(DE(U_0)\right)\right)+D^3E(U_0)\big(DE(U_0),DE(U_0)\big)\big]+\text{h.o.t.}
\end{split}
\end{align}
\normalsize
We now present the error at each stage of the multi-stage algorithm, \cref{eq:ms}, and the conditions required to achieve various orders of accuracy: 
\begin{proposition}
\label{claim:cons}
Let $U_m$ be given in \cref{eq:ms} for $m=0,1,\ldots,M$. The Taylor expansion of $U_m$ at each stage has the same form as \cref{eq:true}, namely: 
\small
\begin{multline}
\label{eq:indte}
    U_m=U_0-\beta_{1,m}k DE(U_0)+\beta_{2,m}k^2D^2E(U_0)DE(U_0)\\-k^3\big[\beta_{3,m}D^2E(U_0)\left(D^2E(U_0)\left(DE(U_0)\right)\right)+\beta_{4,m}D^3E(U_0)\big(DE(U_0),DE(U_0)\big)\big]+\mathcal{O}(k^4)
\end{multline}
\normalsize
where the coefficients obey the following recursive relation 
\begin{align}
\begin{split}
\label{eq:rec}
&\beta_{1,0}=\beta_{2,0}=\beta_{3,0}=\beta_{4,0}=0\\
&\beta_{1,m}=\frac{1}{S_m}\bigg[1+\sum^{m-1}_{i=1}\gamma_{m,i}\beta_{1,i} \bigg]\\
&\beta_{2,m}=\frac{1}{S_m}\bigg[\beta_{1,m}+\sum^{m-1}_{i=1}\gamma_{m,i}\beta_{2,i}\bigg]\\
&\beta_{3,m}=\frac{1}{S_m}\bigg[\beta_{2,m}+\sum^{m-1}_{i=1}\gamma_{m,i}\beta_{3,i}\bigg]\\
&\beta_{4,m}=\frac{1}{S_m}\bigg[\frac{\beta_{1,m}^2}{2}+\sum^{m-1}_{i=1}\gamma_{m,i} \beta_{4,i}\bigg]
\end{split}
\end{align}
with $S_m=\sum^{m-1}_{i=0}\gamma_{m,i}$.
Furthermore, the following conditions for $U_M$ in scheme \cref{eq:ms} are necessary and
sufficient for various orders of accuracy:

\begin{alignat}{3}
\label{eq:cons}
 &\text{\underline{First Order:}}&\quad & \text{\underline{Second Order:}} &\quad &\text{\underline{Third Order:}} \nonumber \\
&\beta_{1,M}=1& &\beta_{1,M}=1& &\beta_{1,M}=1 \nonumber\\
& & &\beta_{2,M}=1/2 & &\beta_{2,M}=1/2\\
& & & & &\beta_{3,M}=1/6 \nonumber\\
& & & & &\beta_{4,M}=1/6 \nonumber
\end{alignat}
\end{proposition}

\begin{proof}
We will now show by induction that the aforementioned consistency formulas, \cref{eq:indte} and  \cref{eq:rec}, hold.\\

\textbf{Stage zero:} $U_0$ trivially satisfies \cref{eq:indte} and \cref{eq:rec}.\\

\textbf{Stage m:} 
Assume \cref{eq:indte} and \cref{eq:rec} up to stage $m-1$. 
First we are going to solve for $U_m-U_0$ in \cref{eq:ms2}:
\small
\begin{align}
\label{eq:stepone}
U_{m}-U_0=&-\frac{k}{S_m}DE(U_m)+\frac{1}{S_m}\sum^{m-1}_{i=0}\gamma_{m,i}U_i-U_0.
\end{align}
\normalsize

 Now Taylor expand $DE(U_m)$ around $U_0$ in \cref{eq:stepone}:
 \normalsize
\begin{align*}
\\U_{m}-U_0=&-\frac{k}{S_m}\bigg[DE(U_0)+D^2E(U_0)(U_m-U_0)+\frac{1}{2}D^3E(U_0)\big(U_m-U_0,U_m-U_0\big)\bigg]\\&+\frac{1}{S_m}\sum^{m-1}_{i=0}\gamma_{m,i}U_i-U_0+\text{h.o.t.}
\end{align*}
\normalsize
Substituting the ansatz $U_0+kA_1+k^2A_2+k^3A_3+\mathcal{O}(k^4)$ for $U_{m}$ and equation \cref{eq:indte} for $U_i$, and retaining up to terms of third order, we have that
\small
\begin{align}
\label{eq:withansatz}
\begin{split}
kA_1&+k^2A_2+k^3A_3=\\&-\frac{k}{S_m}\bigg[1+\sum^{m-1}_{i=1}\gamma_{m,i}\beta_{1,i}\bigg]DE(U_0)+k^2\bigg[\frac{1}{S_m}D^2E(U_0)\big(-A_1+\sum^{m-1}_{i=1}\gamma_{m,i}\beta_{2,i}DE(U_0)\big)\bigg]\\&-k^3\bigg[\frac{1}{S_m}D^2E(U_0)\big(A_2+\sum^{m-1}_{i=1}\gamma_{m,i}\beta_{3,i}D^2E(U_0)DE(U_0)\big)\\
&+\frac{1}{2}\frac{1}{S_m}D^3E(U_0)\big(A_1,A_1\big)+\frac{1}{S_m}\sum^{m-1}_{i=1}\gamma_{m,i}\beta_{4,i}D^3E(U_0)\big(DE(U_0),DE(U_0)\big)\bigg]+\mathcal{O}(k^4)
\end{split}
\end{align}
\normalsize
Solving for $A_1$, $A_2$, $A_3$ by matching terms of the same order in \cref{eq:withansatz}, we arrive at:
\footnotesize
\begin{align*}
\label{eq:solveansatz}
\begin{split}
A_1=&-\frac{1}{S_m}\bigg[1+\sum^{m-1}_{i=1}\gamma_{m,i}\beta_{1,i}\bigg]DE(U_0)\\
A_2=&\frac{1}{S_m}D^2E(U_0)\big(-A_1+\sum^{m-1}_{i=1}\gamma_{m,i}\beta_{2,i}DE(U_0)\big)\\=&\bigg(\frac{1}{S_m^2} \bigg[1+\sum^{m-1}_{i=1}\gamma_{m,i}\beta_{1,i}\bigg]+\frac{1}{S_m}\sum^{m-1}_{i=1}\gamma_{m,i}\beta_{2,i}\bigg)D^2E(U_0)\big(DE(U_0)\big)\\
A_3=&-\frac{1}{S_m}D^2E(U_0)\big(A_2+\sum^{m-1}_{i=1}\gamma_{m,i}\beta_{3,i}D^2E(U_0)DE(U_0)\big)\\
&-\frac{1}{2}\frac{1}{S_m}D^3E(U_0)\big(A_1,A_1\big)-\frac{1}{S_m}\sum^{m-1}_{i=1}\gamma_{m,i}\beta_{4,i}D^3E(U_0)\big(DE(U_0),DE(U_0)\big)\\=-\bigg(&\frac{1}{S_m^3}\bigg[1+\sum^{m-1}_{i=1}\gamma_{m,i}\beta_{1,i} \bigg]+\frac{1}{S_m^2}\sum^{m-1}_{i=1}\gamma_{m,i}\beta_{2,i}\\&+\frac{1}{S_m}\sum^{m-1}_{i=1}\gamma_{m,i}\beta_{3,i}\bigg)D^2E(U_0)\left(D^2E(U_0)\left(DE(U_0)\right)\right)\\-\bigg(&\frac{1}{2}\frac{1}{S_m^3}\bigg[1+\sum^{m-1}_{i=1}\gamma_{m,i}\beta_{1,i} \bigg]^2+\frac{1}{S_m}\sum^{m-1}_{i=1}\gamma_{m,i} \beta_{4,i}\bigg)D^3E(U_0)\big(DE(U_0),DE(U_0)\big)
\end{split}
\end{align*}
\normalsize
completing the induction step. \\

Matching the consistency equations, \cref{eq:indte} and \cref{eq:rec}, at $U_M$ with the one step error \cref{eq:true} gives the conditions on $U_M$ for various orders of accuracy \cref{eq:cons}, completing the proof.
\end{proof}

In the next section, we give examples of $\gamma$'s that satisfy the consistency equations (\cref{claim:cons}) as well as the hypothesis of \cref{claim:ms} concurrently. 

\section{The New Schemes: Examples}
\label{sec:examples}
In this section, we exhibit second and third order examples of scheme \cref{eq:ms} that satisfy concurrently the hypothesis guaranteeing unconditional energy stability (\cref{claim:ms}) and the consistency equations (\cref{claim:cons}) up to second and third order. We found the $\gamma$'s by the following numerical procedure:  we first found a set of $\gamma$'s that satisfied the conditions of \cref{claim:ms}. Then we used the interior point method with the conditions of \cref{claim:ms} as our constraint and an objective function that penalized the mismatch between the current $\beta_{1,M}$ and $\beta_{2,M}$ (and $\beta_{3,M} $ and $\beta_{4,M}$ for third order) and \cref{eq:cons}. After obtaining $\gamma$'s that satisfied the consistency equations up to some small tolerance as well as our constraint,  we sought a nearby algebraic solution to the consistency equations that still satisfied the conditions in \cref{claim:ms}. For some number of stages $M$, it is impossible to satisfy the consistency equation for a given order and the stability conditions. Therefore, we searched for $\gamma$'s that encoded stable algorithms of various orders with different total number of stages and report a set of $\gamma$'s with the lowest number of stages for a given order here. Using this method we were able to find second and third order stable schemes. 
Whether even higher order accuracy (together with stability) can be obtained with this class of schemes will require a more systematic approach to the solvability of the conditions on $\gamma$, and will be the subject of future work.
  \\
\subsection{Second order examples}
It can be shown that there is no unconditionally energy stable second order two-stage method. However, it turns out that three stages are sufficient for unconditional stability: 
\begin{equation}
\label{eq:2ndordergamma}
    \gamma=\left(\begin{array}{ccc}
 \gamma_{1,0} & 0 & 0  \\
  \gamma_{2,0} & \gamma_{2,1} & 0  \\
 \gamma_{3,0} & \gamma_{3,1} & \gamma_{3,2}  \\
\end{array}\right)=\left(
\begin{array}{ccc}
 5 & 0 & 0  \\
 -2 & 6 & 0  \\
 -2 & \frac{3}{14} & \frac{44}{7}  \\
\end{array}\right)\\
\approx \left(
\begin{array}{ccc}
 5.0 & 0 & 0  \\
 -2.0 & 6.0 & 0  \\
 -2.0 & 0.22 & 6.29  \\
\end{array}\right).
\end{equation}
This choice of $\gamma$'s that endows the three-stage method \cref{eq:ms} with unconditional stability and second order accuracy is by no means unique; indeed, here is another that has the additional benefit of having each one of its stages depend only on the previous one and $u_n$:
\begin{equation}
\label{eq:2ndordergamma2}
    \gamma=\left(
\begin{array}{ccc}
 \frac{9}{2} & 0 & 0  \\
 -\frac{11}{6} & \frac{44}{7} & 0  \\
 -\frac{287591}{148306} & 0 & \frac{944163}{148306}  \\
\end{array}\right)\approx \left(
\begin{array}{ccc}
 4.5 & 0 & 0  \\
 -1.83 & 6.29 & 0  \\
 -1.94 & 0 & 6.37  \\
\end{array}\right).\\
\end{equation}
\subsection{Third order examples}
We now exhibit a six stage version of scheme \cref{eq:ms} that concurrently satisfies the conditions for unconditional energy stability (\cref{claim:ms}) as well the consistency equations (\cref{claim:cons})  up to third order:
\begin{align}
\label{eq:3rdordergamma}
&\gamma \approx \left(
\begin{array}{cccccc}
 11.17& 0& 0&0&0&0  \\
 -7.5 & 19.43 & 0&0&0&0  \\
 -1.05& -4.75&13.98&0&0&0 \\
 1.8&0.05&-7.83&13.8&0& 0\\
 6.2&-7.17& -1.33&1.63&11.52&0\\
 -2.83&4.69&2.46&-11.55&6.68&11.95\\
\end{array}\right)
\end{align}
The exact values of the $\gamma$'s above are given in the appendix (\cref{sec:append}); they are all rational numbers but with long fractional representations. Again, we cannot rule out other solutions for $\gamma$, possibly with fewer stages.
\section{The New Schemes: Numerical Tests}
\label{sec:nr}
In this section, we will apply the second order \cref{eq:2ndordergamma} and third \cref{eq:3rdordergamma} order accurate unconditionally stable schemes to a variety of gradient flows. We found \cref{eq:2ndordergamma} before \cref{eq:2ndordergamma2} and therefore ran all our numerical tests with the former. The gradient flows considered span linear and non-linear ordinary and partial differential equations. The corresponding energies include convex and non-convex forms.
Careful numerical convergence studies are presented in each case to verify the anticipated convergence rates of previous sections.
\begin{remark}
Note that equation \cref{eq:ms} can be rewritten using only one quadratic movement limiter term, so a black box implementation for backward Euler \cref{eq:bE}, or equivalently \cref{eq:minmov}, is all that is needed for our method, and is called once per stage. 
\end{remark}
\subsection{Ordinary Differential Equations}
\begin{figure}[ht]
\begin{center}
\includegraphics[width=.45\textwidth]{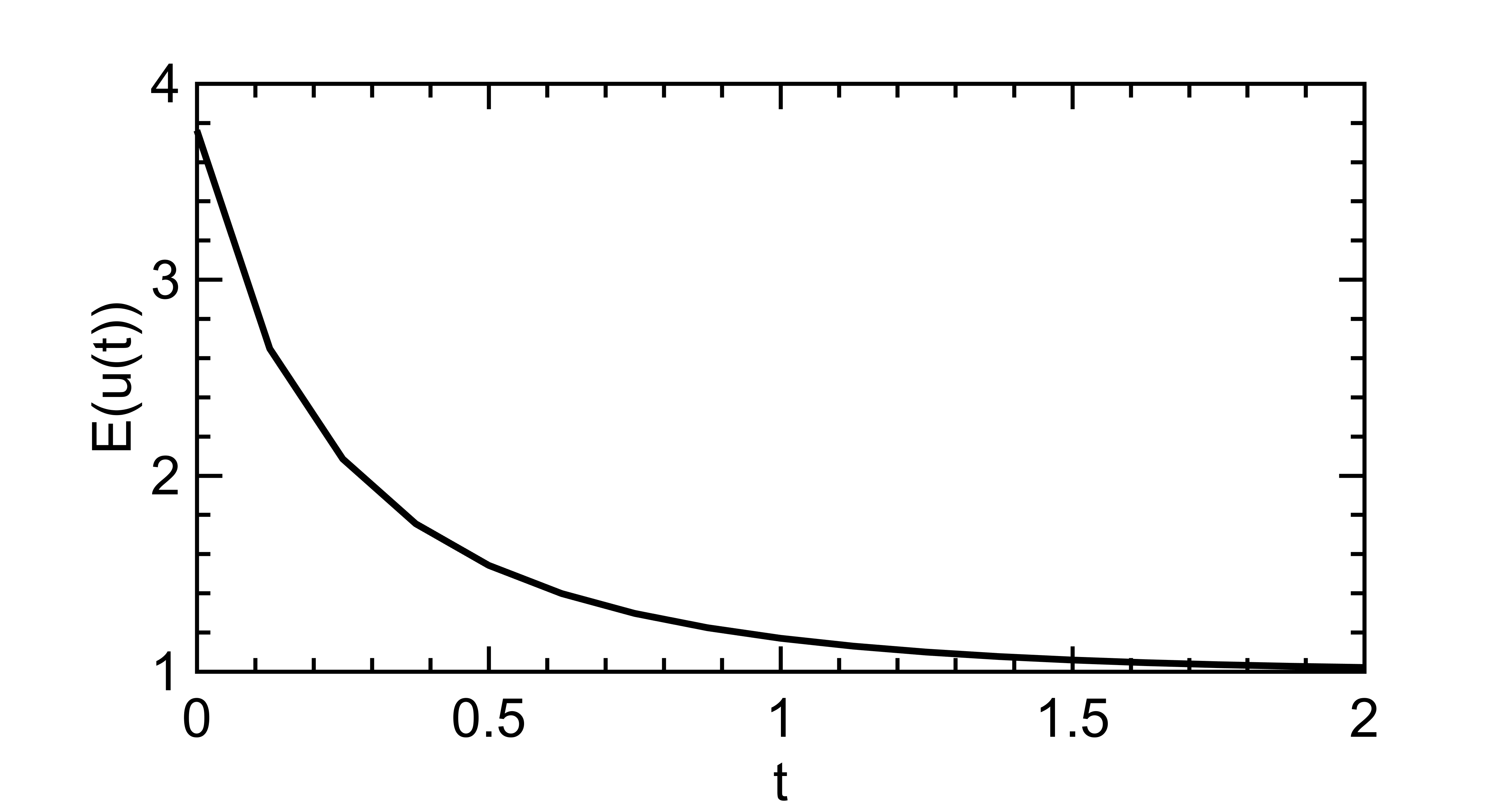}
\caption{\footnotesize The new third six-stage scheme \eqref{eq:ms} \& \eqref{eq:3rdordergamma} with 16 time steps on the ODE on the ODE $u'=-\sinh(u)$ with energy $E(u) = \cosh(u)$}
\label{fig:energyODE}
\end{center}
\end{figure}

First, we turn to the ODE $u' = -\sinh(u)$ with the corresponding energy $E(u)=\cosh(u)$.
With initial condition $u(0)=-2$, the exact solution is $u_*(t)=-2\coth^{-1}(\exp(t)\coth(1))$.
\cref{tab:c3} and \cref{tab:c6} show the error in the solution at time $t=2$ computed by the second order scheme \eqref{eq:ms} \& \eqref{eq:2ndordergamma} and the third order scheme \eqref{eq:ms} \& \eqref{eq:3rdordergamma}, respectively, at various choices of the time step size.
The anticipated order or convergence is clearly observed for both schemes. \cref{fig:energyODE} shows the energy at every time step for the third order method with 16 time steps. As expected, the energy decreases at every time step. There is little visual difference between \cref{fig:energyODE}, the plot of the second order method with 16 time steps and the plot of the exact energy.


\begin{table}[ht]
\begin{center}
\begin{tabular}{|c|c|c|c|c|c|}
\hline
Number of& & & & &\\time steps&$2^{4}$&$2^{5}$&$2^{6}$&$2^{7}$&$2^{8}$\\
\hline
Error at $t=2$& 5.25e-04&1.31e-04&3.27e-05&8.18e-06&2.05e-06\\
\hline
Order&-&2.00&2.00&2.00&2.00\\
\hline
\end{tabular}
\caption{\footnotesize The new second order accurate, unconditionally stable, three-stage scheme \eqref{eq:ms} \& \eqref{eq:2ndordergamma} on the ODE $u'=-\sinh(u)$ with energy $E(u) = \cosh(u)$.}
\label{tab:c3}
\end{center}
\end{table}

\begin{table}[ht]
\begin{center}
\begin{tabular}{|c|c|c|c|c|c|}
\hline
Number of& & & & &\\time steps&$2^{4}$&$2^{5}$&$2^{6}$&$2^{7}$&$2^{8}$\\
\hline
Error at $t=2$& 1.19e-05&1.48e-06&1.85e-07&2.30e-08&2.88e-09\\
\hline
Order&-&3.00&3.00&3.00&3.00\\
\hline
\end{tabular}
\caption{\footnotesize The new third order accurate, unconditionally stable, six-stage scheme \eqref{eq:ms} \& \eqref{eq:3rdordergamma} on the ODE $u'=-\sinh(u)$ with energy $E(u) = \cosh(u)$.}
\label{tab:c6}
\end{center}
\end{table}

\begin{figure}[ht]
\begin{center}
\includegraphics[width=.45\textwidth]{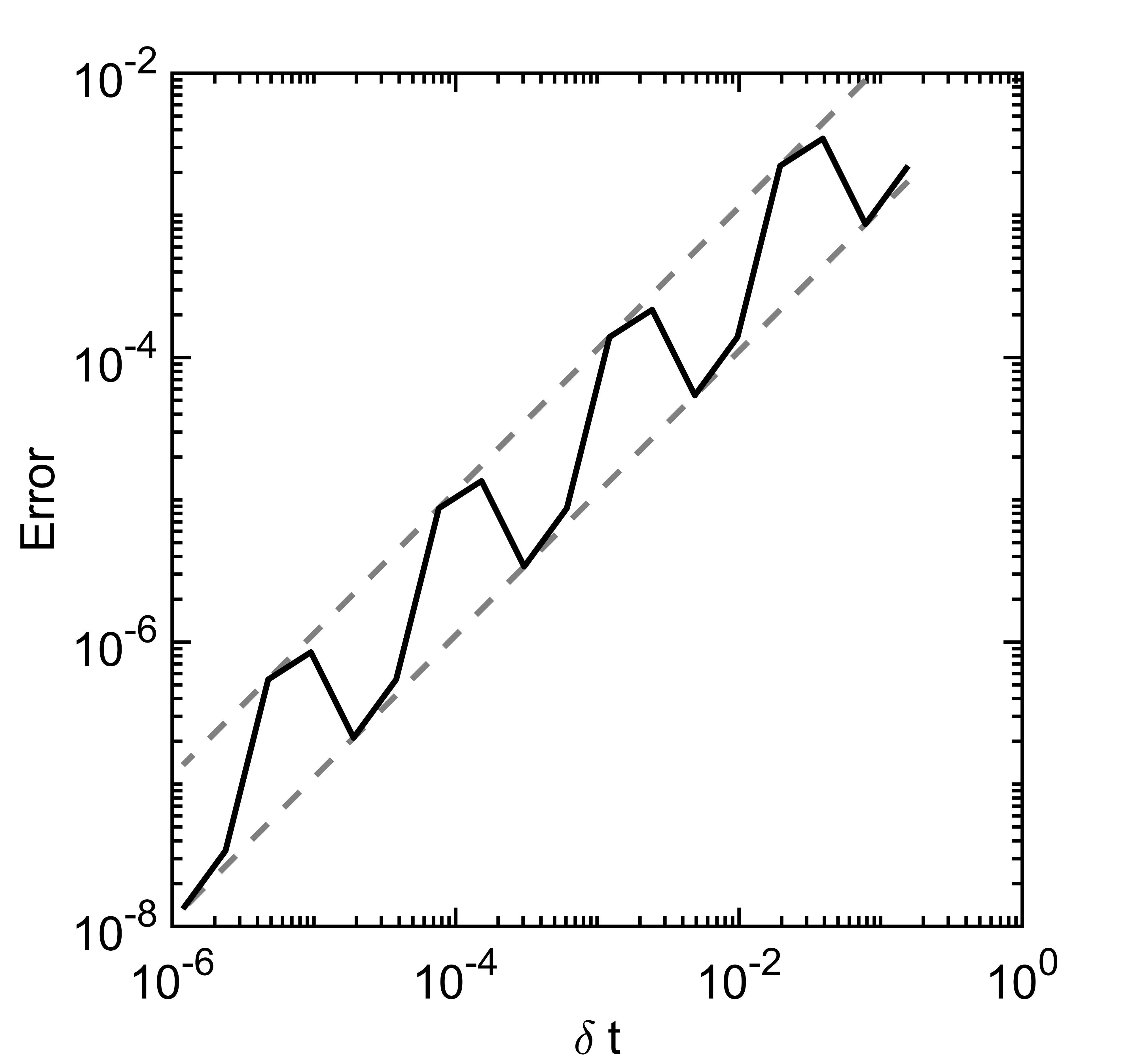}
\includegraphics[width=.45\textwidth]{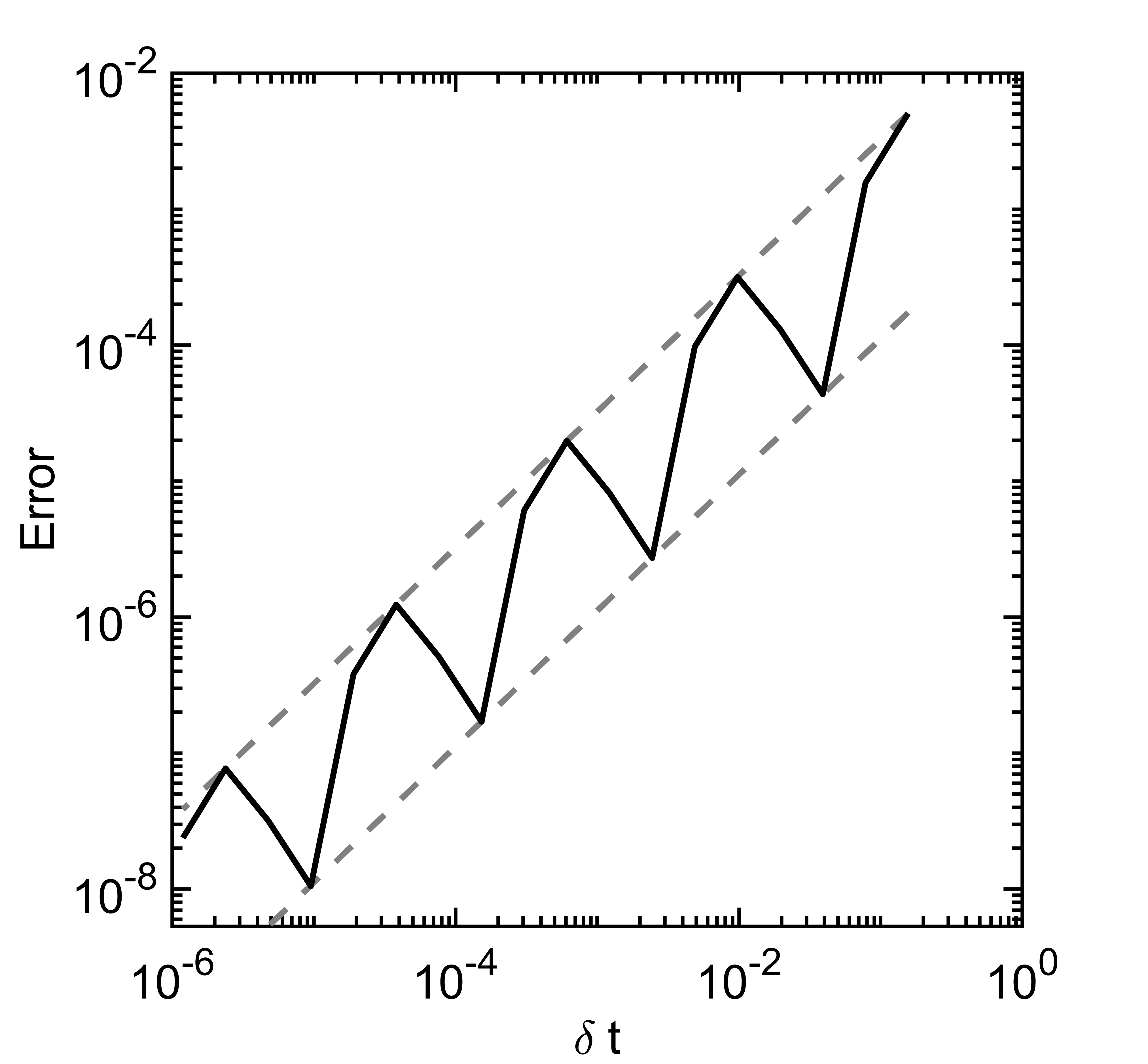}
\caption{\footnotesize The new second order accurate, unconditionally stable, three-stage scheme \eqref{eq:ms} \& \eqref{eq:3rdordergamma} (right) and the new third six-stage scheme \eqref{eq:ms} \& \eqref{eq:3rdordergamma} (left) on the ODE induced by gradient flow on non smooth energy \cref{eq:nonsmooth}}
\label{fig:nonsmooth}
\end{center}
\end{figure}

\begin{figure}[ht]
\begin{center}
\includegraphics[width=.45\textwidth]{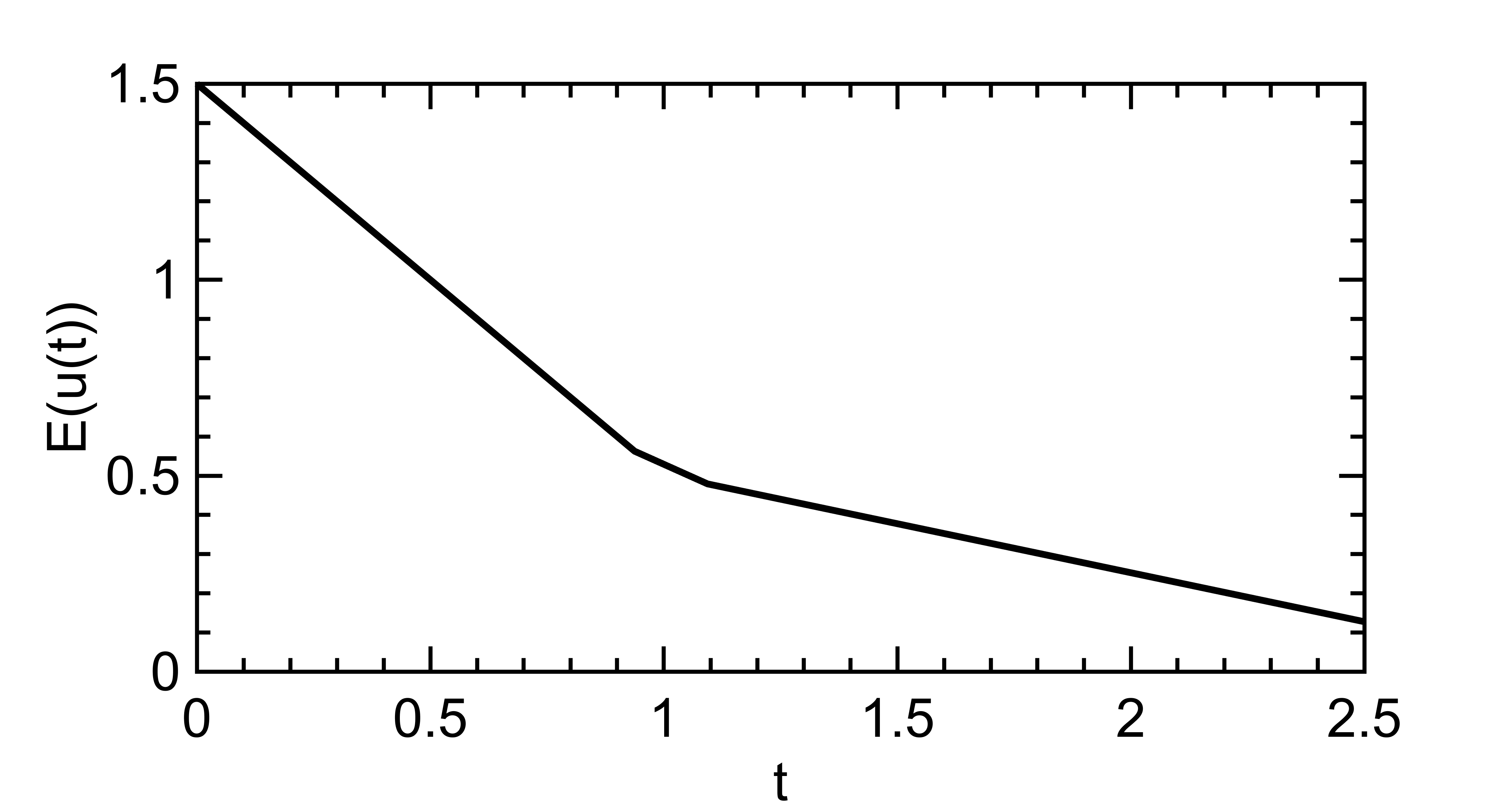}
\caption{\footnotesize The new third order six-stage scheme \eqref{eq:ms} \& \eqref{eq:3rdordergamma} with 16 time steps on the ODE induced by gradient flow on non smooth energy \cref{eq:nonsmooth}}
\label{fig:nonsmoothenergy}
\end{center}
\end{figure}

We next turn to an ODE with the non-smooth energy
\begin{equation}
\label{eq:nonsmooth}
  E(u) =
  \begin{cases}
                                   \frac{1}{2}|u| & \text{if $|u|<1$} \\
                                   |u-1|+\frac{1}{2} & \text{if $|u|\geq 1$}.
  \end{cases}
\end{equation}
Since the energy is non-smooth we do not expect higher order convergence. As shown in \cref{fig:nonsmooth}, the second \eqref{eq:ms} \& \eqref{eq:2ndordergamma} and third order scheme \eqref{eq:ms} \& \eqref{eq:3rdordergamma} obtain first order convergence on average. Notwithstanding, the energy decreases at every time step as shown in \cref{fig:nonsmoothenergy} for the third order method with 16 time steps.

\subsection{Partial Differential Equations}
\label{subsec:PDE}
\begin{figure}[ht]
\begin{center}
\includegraphics[width=.45\textwidth]{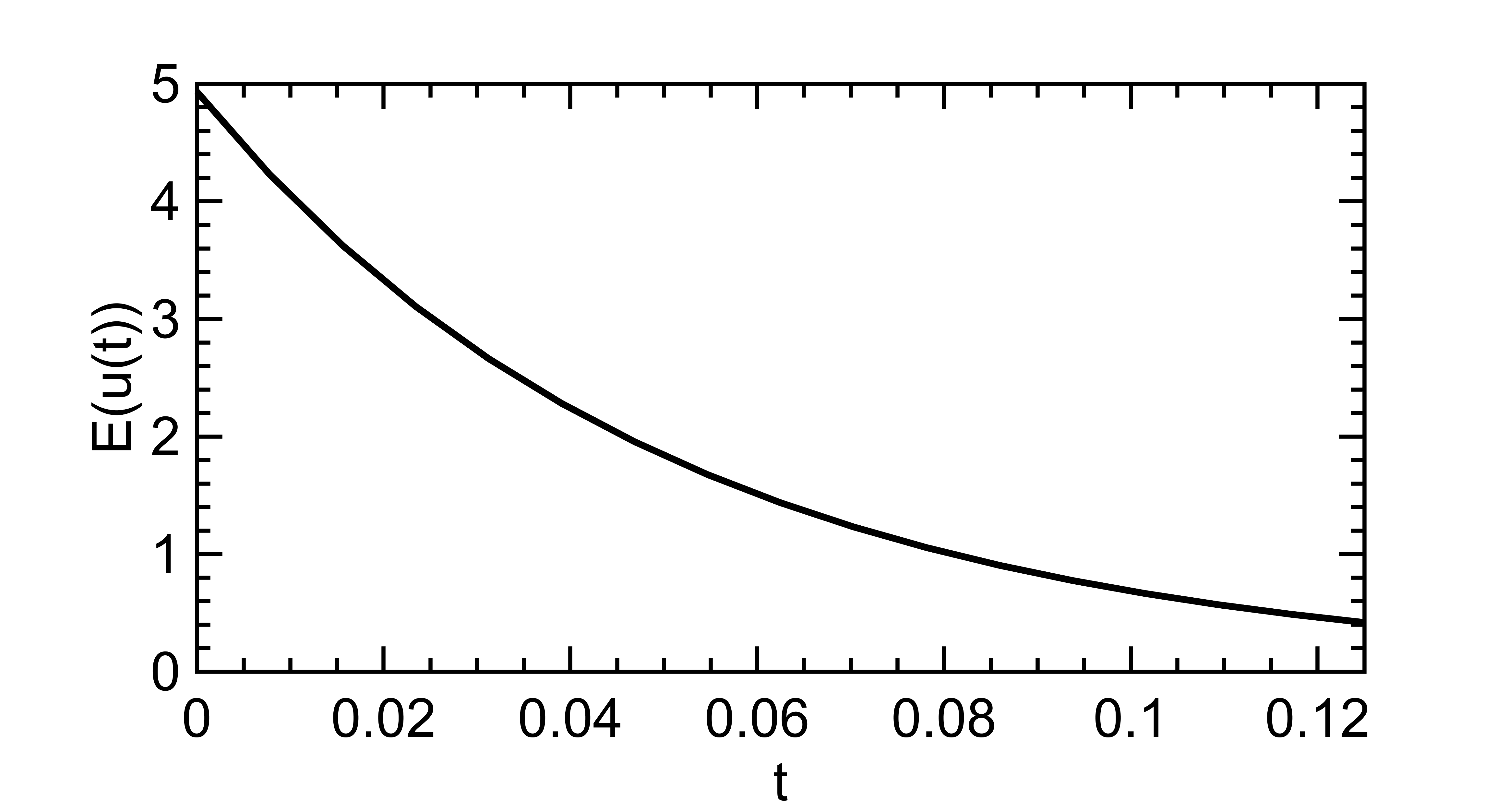}
\caption{\footnotesize The new third six-stage scheme \eqref{eq:ms} \& \eqref{eq:3rdordergamma} with 16 time steps on the one-dimensional heat equation $u_t = u_{xx}$}
\label{fig:energyheat}
\end{center}
\end{figure}

For PDEs, we start with a preliminary test on the one dimensional heat equation $u_t = u_{xx}$ on $x\in[-1,1]$ subject to periodic boundary conditions with initial data $u(x,0)=\sin(\pi x)$.
This is gradient flow with respect to the $L^2$ inner product for the energy $E(u) = \frac{1}{2} \int u_x^2 \, dx$.
The exact solution is $u_*(x,t)=\sin(\pi x)\exp(-\pi^2 t)$.
The spatial domain $[-1,1]$. For this example as well as the other PDEs in this section, we choose the discretization of the Laplacian and number of spatial points so that the contribution to the error from the spatial discretization is negligible.  
 \cref{tab:h3} and \cref{tab:h6} show the $L^2$ error in the approximate solution at $t=\frac{1}{8}$, computed via the second order accurate scheme \eqref{eq:ms} \& \eqref{eq:2ndordergamma}, and the third order accurate scheme \eqref{eq:ms} \& \eqref{eq:3rdordergamma}, respectively.
 \cref{fig:energyheat} shows the energy at every time step for the third order method with 16 time steps. We see that the energy decreases at every time step.

\begin{table}[ht]
\begin{center}
\begin{tabular}{|c|c|c|c|c|c|c|}
\hline
Number of& & & & & & \\time steps& $2^2$&$2^{3}$&$2^{4}$&$2^{5}$&$2^{6}$&$2^{7}$\\
\hline
L2& 1.09e-03&2.66e-04&6.59e-05&1.64e-05&4.09e-06 &1.02e-06\\
\hline
Order&-&2.03&2.01&2.01&2.00&2.00\\
\hline
\end{tabular}
\caption{\footnotesize The new second order accurate, unconditionally stable, three-stage scheme \eqref{eq:ms} \& \eqref{eq:2ndordergamma} on the one-dimensional heat equation $u_t = u_{xx}$.}
\label{tab:h3}
\end{center}
\end{table}

\begin{table}[ht]
\begin{center}
\begin{tabular}{|c|c|c|c|c|c|c|}
\hline
Number of& & & & & &\\time steps& $2^2$&$2^{3}$&$2^{4}$&$2^{5}$&$2^{6}$&$2^{7}$\\
\hline
L2&2.30e-05& 2.75e-06&3.36e-07&4.16e-06&5.17e-09&6.37e-10\\
\hline
Order&-&3.06&3.03&3.02&3.01&3.02\\
\hline
\end{tabular}
\caption{\footnotesize The new third order accurate, unconditionally stable, six-stage scheme \eqref{eq:ms} \& \eqref{eq:3rdordergamma} on the one-dimensional heat equation $u_t = u_{xx}$.}
\label{tab:h6}
\end{center}
\end{table}

\begin{figure}[ht]
\begin{center}
\includegraphics[width=.45\textwidth]{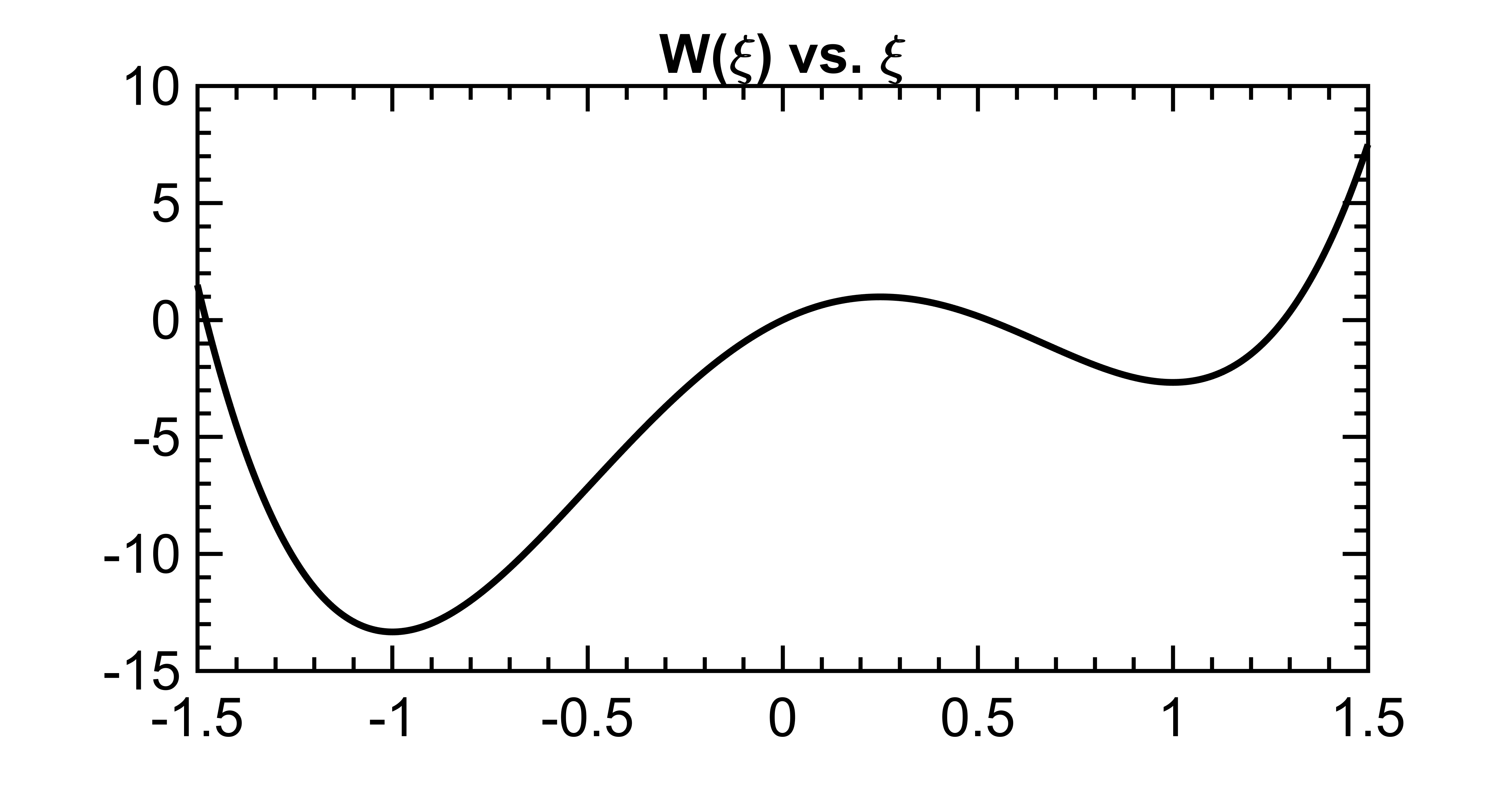}
\includegraphics[width=.45\textwidth]{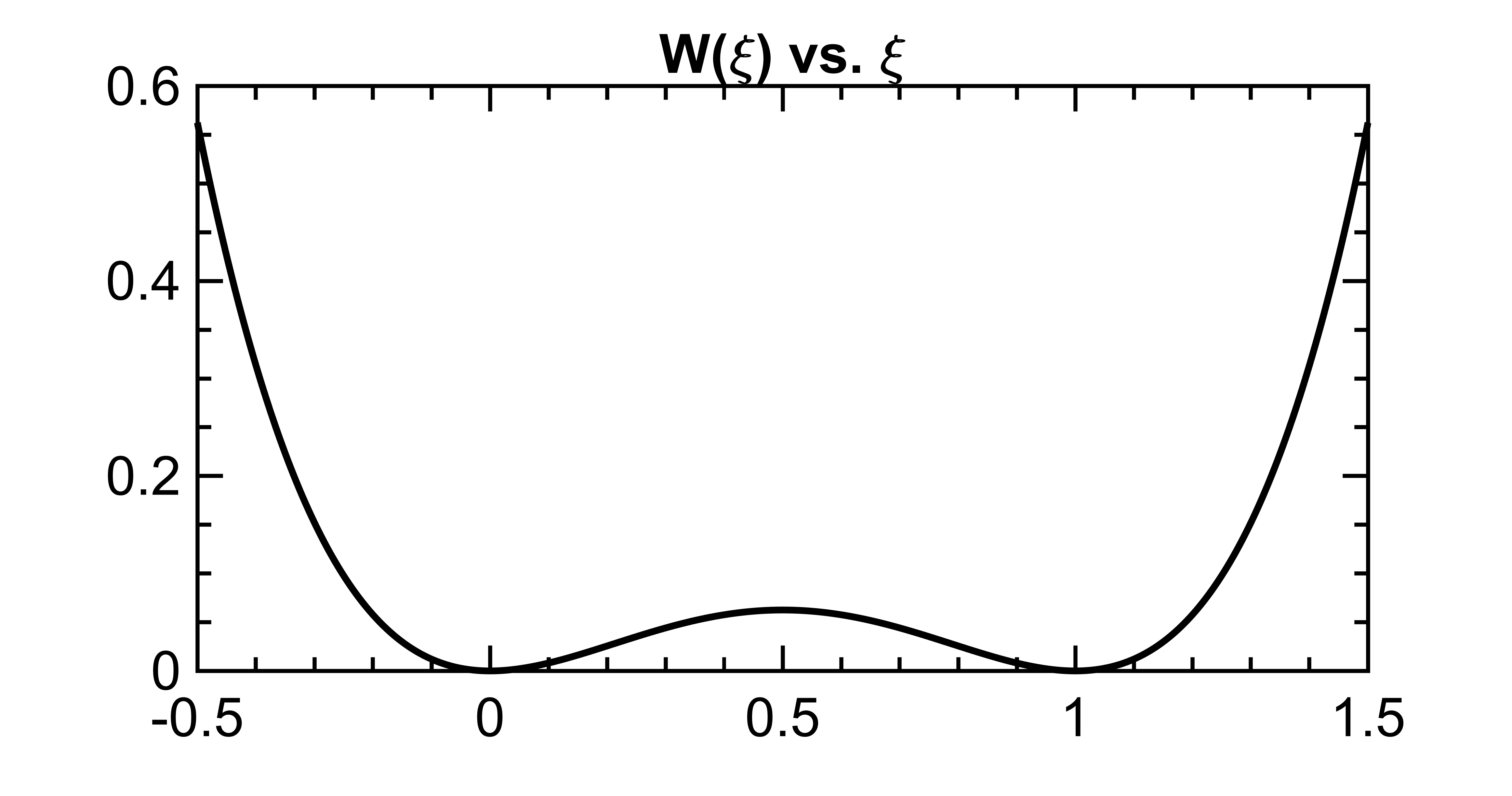}
\caption{\footnotesize The double well potentials used in the Allen-Cahn (\ref{eq:ac}) and Cahn-Hilliard (\ref{eq:ch}) equations: One with unequal and the other with equal depth wells.}
\label{fig:W}
\end{center}
\end{figure}

\begin{figure}[ht]
  \begin{center}
\includegraphics[width=.5\textwidth]{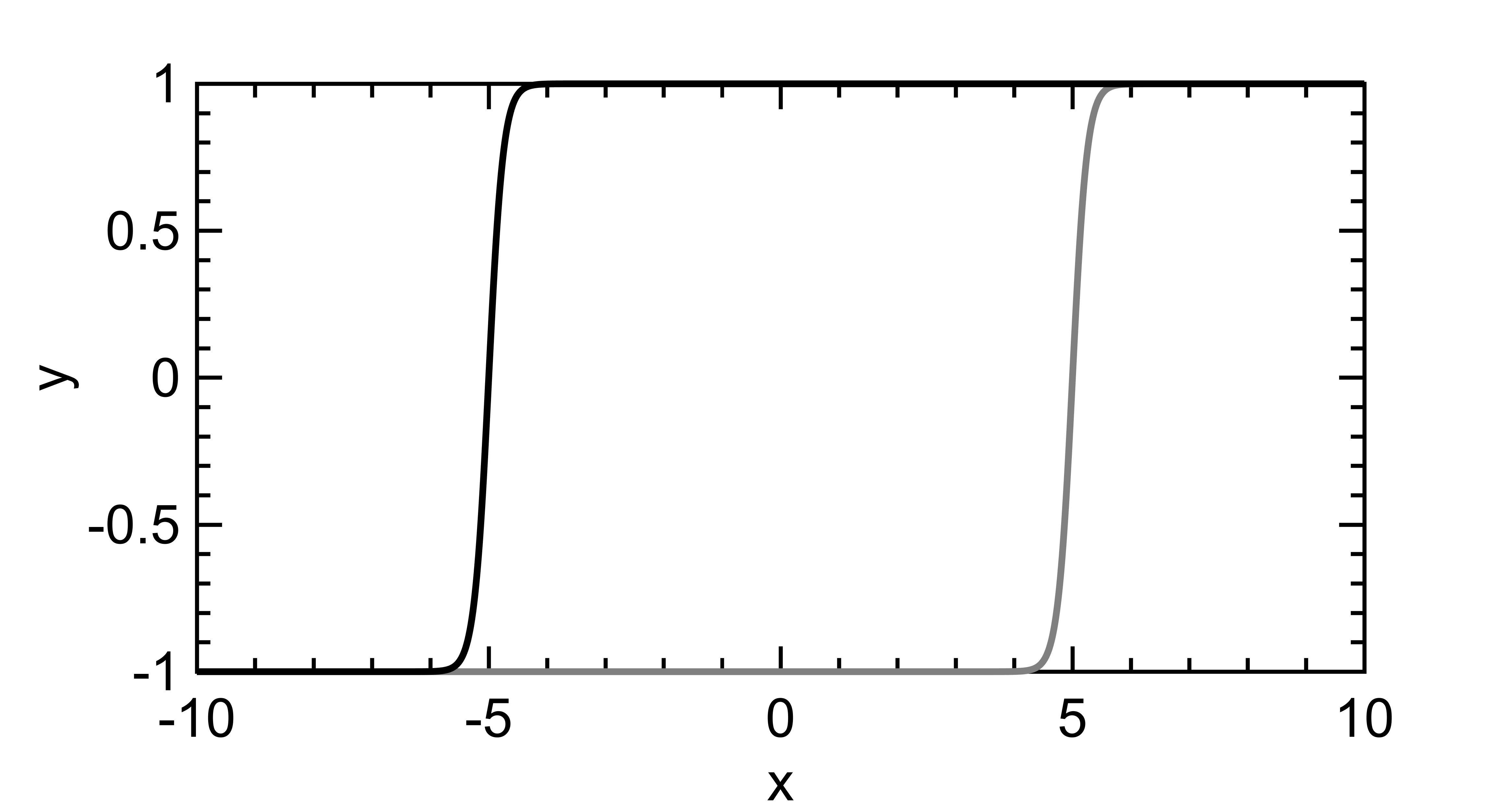} %
\end{center}
\caption{\footnotesize The initial condition (black) and the solution at final time (gray) in the numerical convergence study on the 1D Allen-Cahn equation (\ref{eq:ac}) with a potential that has unequal depth wells.}
\label{fig:1dac}
\end{figure}

We now turn to less trivial examples, starting with the Allen-Cahn equation
\begin{equation}
\label{eq:ac}
u_t = \Delta u - W'(u)
\end{equation}
where $W:\mathbb{R}\to\mathbb{R}$ is a double-well potential.
This is gradient flow for the energy
\begin{equation}
\label{eq:acE}
E(u) = \int \frac{1}{2} |\nabla u|^2 + W(u) \, dx
\end{equation}
with respect to the $L^2$ inner product.

First, we consider equation (\ref{eq:ac}) in one space dimension, with the potential $W(u) = 8u-16u^2-\frac{8}{3}u^3+8u^4$.
This is a double well potential with unequal depth wells; see Figure \ref{fig:W}.
In this case, equation (\ref{eq:ac}) is well-known to possess traveling wave solutions on $x\in\mathbb{R}$, see \cref{fig:1dac}.
We choose the initial condition $u(x,0) = \tanh(4x + 20)$; the exact solution is then $u_*(x,t) = \tanh(4x + 20 - 8t)$.
The computational domain is $x\in[-10,10]$.
We approximate the solution on $\mathbb{R}$ by using the Dirichlet boundary conditions $u(\pm 10,t) = \pm 1$:
The domain size is large enough that the mismatch in boundary conditions do not substantially contribute to the error in the approximate solution over the time interval $t\in[0,5]$.
\cref{tab:twms3} and \cref{tab:twms6} tabulate the error in the computed solution at time $t=5$ for our two new schemes. 

\begin{table}
\begin{center}
\begin{tabular}{|c|c|c|c|c|c|c|c|}
\hline
Number of& & & & & &\\time steps&$2^7$&$2^8$&$2^9$&$2^{10}$&$2^{11}$&$2^{12}$\\
\hline
L2&5.14e-02&1.26e-02& 3.13e-03&7.79e-04&1.94e-04&4.86e-05\\
\hline
Order&-&2.02&2.01&2.01 &2.00 &2.00 \\
\hline
\end{tabular}
\caption{\footnotesize The new second order accurate, unconditionally stable, three-stage scheme \eqref{eq:ms} \& \eqref{eq:2ndordergamma} on the one-dimensional Allen-Cahn equation (\ref{eq:ac}) with a traveling wave solution.}
\label{tab:twms3}
\end{center}
\end{table}

\begin{table}[ht]
\begin{center}
\begin{tabular}{|c|c|c|c|c|c|c|}
\hline
Number of& & & & & &\\time steps&$2^7$&$2^8$&$2^9$&$2^{10}$&$2^{11}$&$2^{12}$\\
\hline
L2&9.06e-04&9.97e-05& 1.20e-05&1.48e-06&1.85e-07&2.37e-08\\
\hline
Order&-&3.18&3.06&3.02 &3.00 &2.97\\
\hline
\end{tabular}
\caption{\footnotesize The new third order accurate, unconditionally stable, six-stage scheme \eqref{eq:ms} \& \eqref{eq:3rdordergamma} on the one-dimensional Allen-Cahn equation (\ref{eq:ac}) with a traveling wave solution.}
\label{tab:twms6}
\end{center}
\end{table}

\begin{figure}[ht]
  \begin{center}
\includegraphics[width=.45\textwidth]{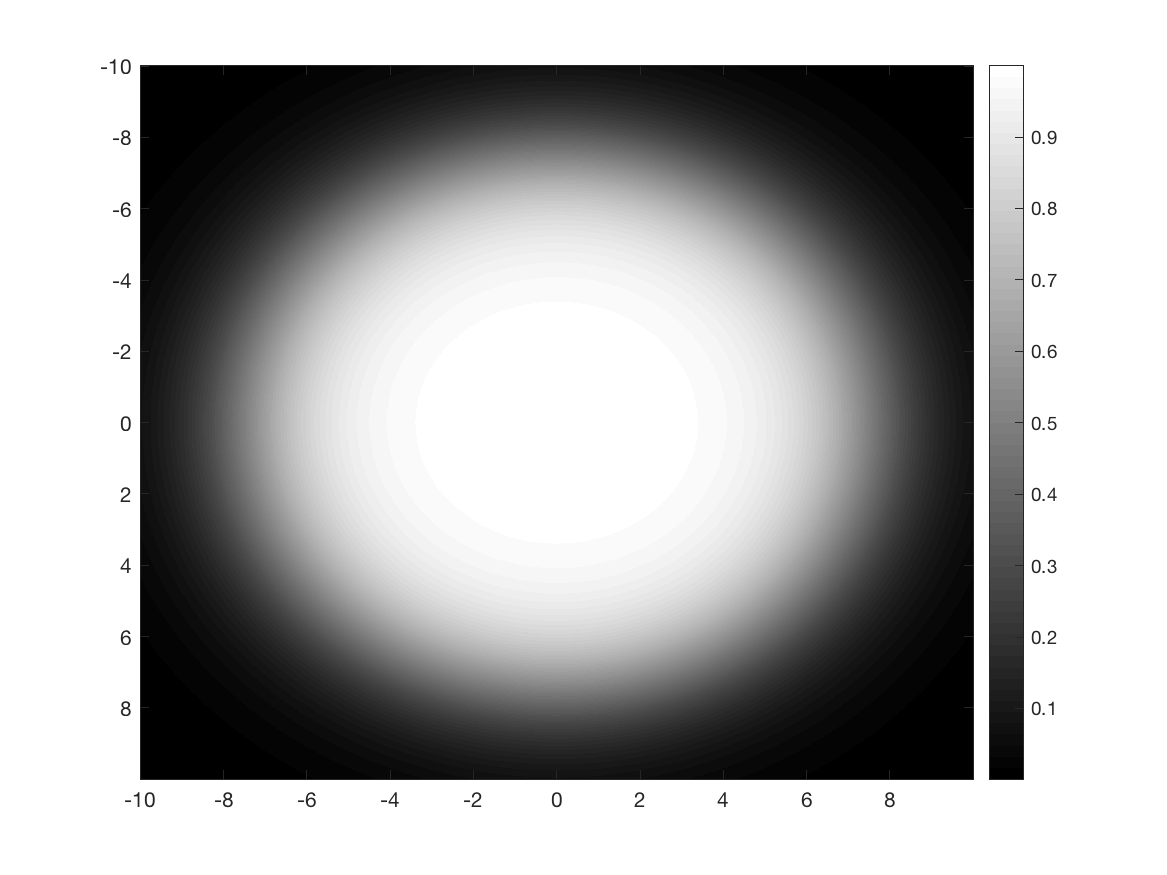}
\includegraphics[width=.45\textwidth]{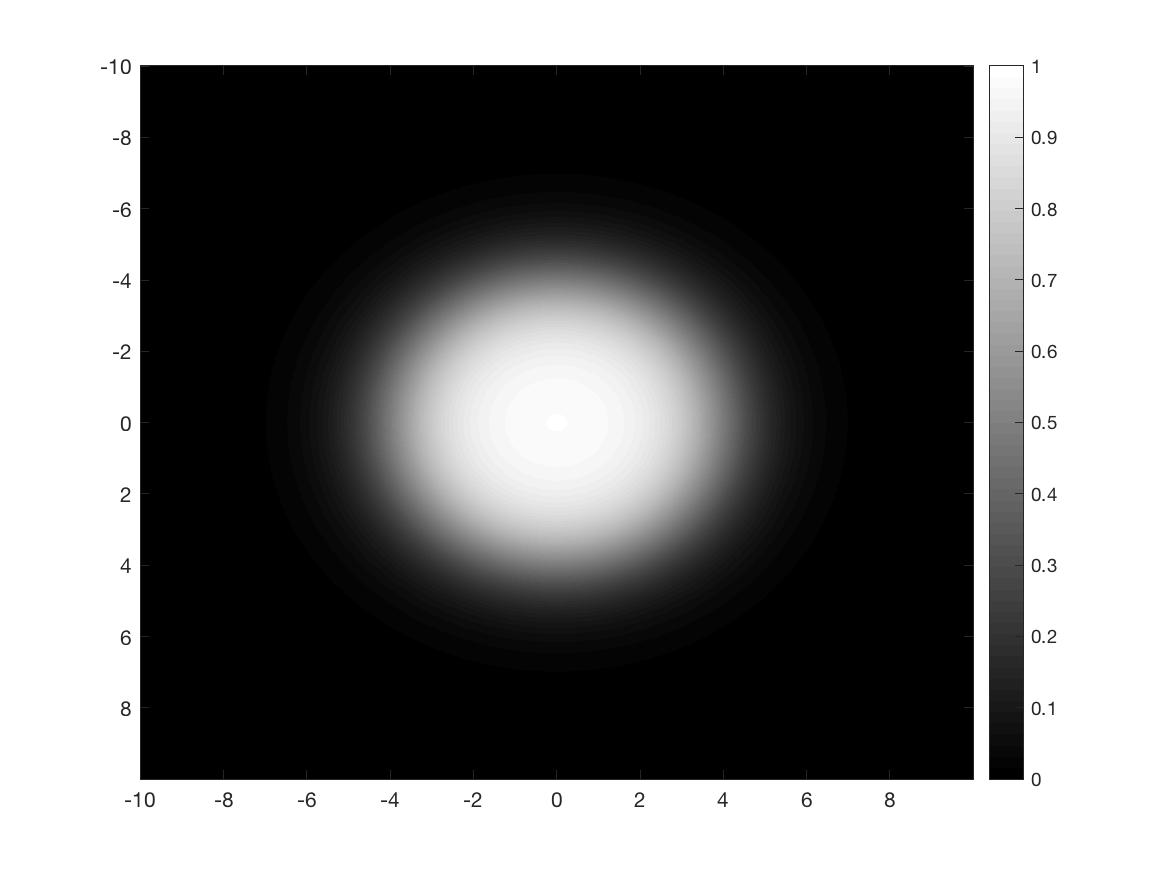}
\end{center}
\caption{\footnotesize Initial condition and the solution at final time for the 2D Allen-Cahn equation with a potential that has equal depth wells.}
\label{fig:2dac}
\end{figure}

Next, we consider the Allen-Cahn equation (\ref{eq:ac}) in two space dimensions, with the potential $W(u) = u^2(1-u)^2$ that has equal depth wells; see Figure \ref{fig:W}.
We take the initial condition $u(x,y,0)=\frac{1}{1+\exp[-(7.5-\sqrt{x^2+y^2})]}$ on the domain $x\in [-10,10]^2$, and impose periodic boundary conditions. We run the system to find $u$ at $t=20$, (\cref{fig:2dac} shows $u$ at $t=0$ and $t=20$).
As a proxy for the exact solution of the equation with this initial data, we compute a very highly accurate numerical approximation $u_*(x,y,t)$ via the following second order accurate in time, semi-implicit, multi-step scheme ~\cite{chen1998applications} on an extremely fine spatial grid and take very small time steps:
\begin{equation*}
    \frac{3}{2}u^{n+1}-2u^{n}+\frac{1}{2}u^{n-1}=k\Delta u^{n+1}-k(2W'(u^{n})-W'(u^{n-1})).
\end{equation*}
\cref{tab:acms32d} and \cref{tab:acms62d} show the errors in and convergence rates for the approximate solutions computed by our new multi-stage schemes.

\begin{table}
\begin{center}
\begin{tabular}{|c|c|c|c|c|c|c|}
\hline
Number of& & & & & &\\time steps&$2^5$&$2^6$&$2^7$&$2^{8}$&$2^{9}$&$2^{10}$\\
\hline
L2&3.43e-03&8.73e-04& 2.21e-04&5.55e-05&1.39e-05&3.49e-06\\
\hline
Order&-&1.98&1.98&1.99 &1.99&2.00 \\
\hline
\end{tabular}
\caption{\footnotesize The new second order accurate, unconditionally stable, three-stage scheme \eqref{eq:ms} \& \eqref{eq:2ndordergamma} on the two-dimensional Allen-Cahn equation (\ref{eq:ac}) with a potential that has equal depth wells.}
\label{tab:acms32d}
\end{center}
\end{table}

\begin{table}[ht]
\begin{center}
\begin{tabular}{|c|c|c|c|c|c|c|}
\hline
Number of& & & & & &\\time steps&$2^3$&$2^4$&$2^5$&$2^{6}$&$2^{7}$&$2^8$\\
\hline
L2&4.60e-03&5.41e-04& 6.44e-05&7.98e-06&1.02e-06&1.33e-07\\
\hline
Order&-&3.09&3.07&3.01 &2.97&2.94 \\
\hline
\end{tabular}
\caption{\footnotesize The new third order accurate, unconditionally stable, six-stage scheme \eqref{eq:ms} \& \eqref{eq:3rdordergamma} on the two-dimensional Allen-Cahn equation (\ref{eq:ac}) with a potential that has equal depth wells.}
\label{tab:acms62d}
\end{center}
\end{table}

\begin{figure}[ht]
  \begin{center}
\includegraphics[width=.45\textwidth]{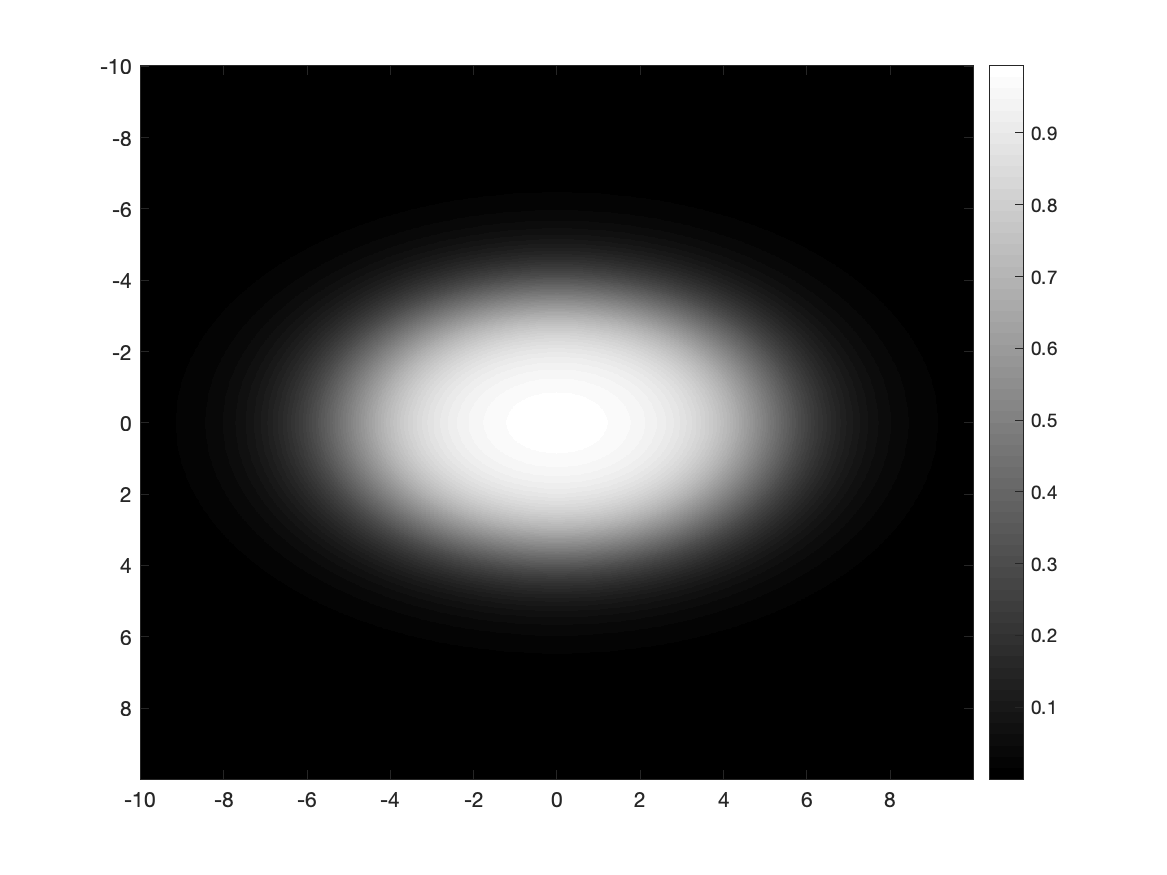}
\includegraphics[width=.45\textwidth]{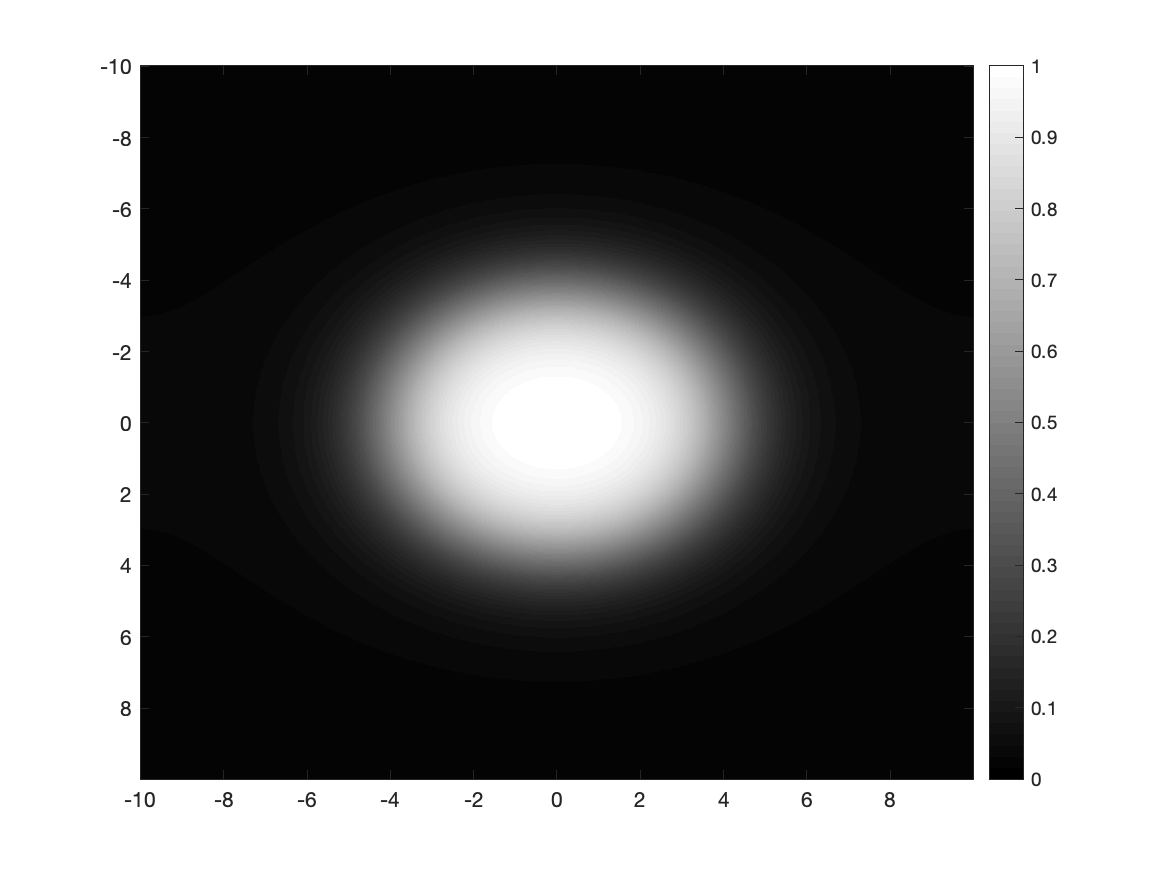}
\end{center}
\caption{\footnotesize Initial condition and the solution at final time for the 2D Cahn-Hilliard equation with a potential that has equal depth wells.}
\label{fig:2dch}
\end{figure}

As a final example, we consider the Cahn-Hilliard equation
\begin{equation}
\label{eq:ch}
u_t = -\Delta \big( \Delta u - W'(u) \big)
\end{equation}
where we take $W$ to be the double well potential $W(u) = u^2(1-u)^2$ with equal depth wells and impose periodic boundary conditions.
This flow is also gradient descent for energy \cref{eq:acE}, but with respect to the $H^{-1}$ inner product:
\begin{equation*}
\langle u \,, v \, \rangle = \int u \Delta^{-1} v \, dx.
\end{equation*}
Starting from the initial condition $u(x,y,0)=\frac{1}{1+\exp[-(5-\sqrt{x^2+2y^2})]}$ and running the system until $t=20$ (see \cref{fig:2dch}). We computed a proxy for the ``exact'' solution once again using the second order accurate, semi-implicit multi-step scheme from \cite{chen1998applications} \cite{shen2010numerical}:
\begin{equation*}
    \frac{3}{2}u^{n+1}-2u^{n}+\frac{1}{2}u^{n-1}=-k\Delta[\Delta u^{n+1}-k(2W'(u^{n})-W'(u^{n-1}))]
\end{equation*}
where the spatial and temporal resolution was taken to be high to ensure the errors are low.
\cref{tab:chms32d} and \cref{tab:chms62d} show the errors in and convergence rates for the approximate solutions computed by our new multi-stage schemes.

\begin{table}
\begin{center}
\begin{tabular}{|c|c|c|c|c|c|}
\hline
Number of& & & & &\\time steps&$2^2$&$2^3$&$2^4$&$2^{5}$&$2^{6}$\\
\hline
L2&1.16e-03&2.62e-04& 6.41e-05&1.64e-05&4.22e-06\\
\hline
Order&-&2.15&2.03&1.97&1.95 \\
\hline
\end{tabular}
\caption{\footnotesize The new second order accurate, unconditionally stable, three-stage scheme \eqref{eq:ms} \& \eqref{eq:2ndordergamma} on the two-dimensional Cahn-Hilliard equation (\ref{eq:ch}) with a potential that has equal depth wells.}
\label{tab:chms32d}
\end{center}
\end{table}

\begin{table}[ht]
\begin{center}
\begin{tabular}{|c|c|c|c|c|}
\hline
Number of& & & & \\time steps&$2^2$&$2^3$&$2^4$&$2^{5}$\\
\hline
L2&2.20e-04 &4.12e-05& 6.73e-06&1.05e-06\\
\hline
Order&-&2.42&2.62&2.67  \\
\hline
\end{tabular}
\caption{\footnotesize The new third order accurate, unconditionally stable, six-stage scheme \eqref{eq:ms} \& \eqref{eq:2ndordergamma} on the two-dimensional Cahn-Hilliard equation (\ref{eq:ch}) with a potential that has equal depth wells.}
\label{tab:chms62d}
\end{center}
\end{table}

\begin{remark}
\label{remark:thresh}
As further evidence of the generality and flexibility of the new schemes introduced in this paper, we note that they can also be used to jack up the order of accuracy in time of less conventional numerical algorithms such as {\it threshold dynamics} \cite{mbo92,mbo94}.
Also known as diffusion generated motion, threshold dynamics is an unconditionally stable algorithm for simulating the motion of interfaces by mean curvature, merely by alternating the two simple steps of convolution and thresholding.
It was given a variational formulation in \cite{eo} that exhibits it as carrying out an approximate minimizing movements procedure at every time step.
Although the stability calculation of Section \ref{sec:gf} applies verbatim, the consistency calculations of Section \ref{sec:ex} have to be redone.
This is because (a) motion by mean curvature, although formally a gradient flow on perimeter, does not quite fit the classical formulation (\ref{eq:de}), and (b) the variational formulation in \cite{eo} shows that threshold dynamics carries out minimizing movements for {\it approximately} the right energy with respect to {\it approximately} the right metric: these additional errors have to be taken into account.
Due to the substantial modifications to the consistency calculation required, extension of the new schemes to enhancing the order of accuracy of threshold dynamics will be taken up in a subsequent, separate paper \cite{alex2019second}.

\end{remark}

 \section{Conclusion}
We presented a class of unconditionally stable, high order in time schemes for gradient flows.
The new schemes can be thought of as a variational analogue of Richardson extrapolation: they enable jacking up the order of accuracy of standard backward Euler method, while maintaining its unconditional stability, at the expense of taking multiple backward Euler time substeps per full time step.
What results is a universal method to jack up the accuracy to at least third order in time whenever a blackbox implementation of the standard backward Euler scheme is available, while increasing overall complexity by only a constant factor.
We demonstrated the method and its advertised accuracy on a number of linear and nonlinear ODEs and PDEs.

Whether this class of schemes can be used to achieve arbitrarily high (i.e. $\geq 4$) order in time accuracy will be the topic of a future investigation.
 \section{Appendix}
 \label{sec:append}
 We record here the exact values for the coefficients $\gamma$ in the six-stage, third order accurate scheme introduced in Section \ref{sec:examples}.
They are rational numbers, but the irreducible fraction representation of some of them are quite long, and were therefore approximated above.
With the universal, exact values given below, we can rigorously state that the new scheme introduced in this paper can be used to jack up the order of accuracy in time of any backward Euler scheme (\ref{eq:bE}) for gradient flows (\ref{eq:de}) to third order while maintaining unconditional energy stability.
 
The matrix of values is:
 \begin{align*}
&\gamma=\left(
\begin{array}{cccccc}
 \frac{67}{6}& 0& 0&0&0&0  \\
 -\frac{15}{2} & \frac{136}{7} & 0&0&0&0  \\
 -\frac{21}{20}& -\frac{19}{4}&\frac{587}{42}&0&0&0 \\
 \frac{9}{5}& \frac{1}{21}&-\frac{47}{6}&\frac{69}{5}&0& 0\\
 \frac{31}{5}&-\frac{43}{6}& -\frac{4}{3}&\frac{13}{8}&\frac{242}{21}&0\\
 -\frac{17}{6}&\frac{75}{16}&\gamma_{6,2}&\gamma_{6,3}&\gamma_{6,4}&\gamma_{6,5}\\
\end{array}\right)
\end{align*}
where
\begin{align*}
&\gamma_{6,2}=-\frac{96877768305591883216465260738322381995331343806720345}{39417514787340924198452679823989476266149744556295712}\\
&\gamma_{6,3}=-\frac{910677500903250179715877776918800480038125970511673389}{78835029574681848396905359647978952532299489112591424}\\
&\gamma_{6,4}=\frac{2985416726242784122189204876225493950575679989899779}{446910598495928845787445349478338733176300958688160}\\
&\gamma_{6,5}=\frac{523180952458721016795516949849623944572931703979520653}{43797238652601026887169644248877195851277493951439680}
\end{align*}
It can be checked that these $\gamma$'s satisfy the inequalities in the hypothesis of \cref{claim:ms} for stability, and the consistency equations in \cref{claim:cons} for third order exactly. Code for doing so can be found at \url{https://github.com/AZaitzeff/gradientflow}.
\bibliographystyle{siamplain}
\bibliography{references}

\begin{thebibliography}{10}

\bibitem{butcher2016numerical}
{\sc J.~Butcher}, {\em Numerical Methods for Ordinary Differential Equations},
  Wiley, 2016.

\bibitem{butcher1975}
{\sc J.~C. Butcher}, {\em A stability property of implicit {R}unge-{K}utta
  methods}, {BIT}, 15 (1975), pp.~358--361.

\bibitem{chen1998applications}
{\sc L.~Q. Chen and J.~Shen}, {\em Applications of semi-implicit
  {Fourier}-spectral method to phase field equations}, Computer Physics
  Communications, 108 (1998), pp.~147--158.

\bibitem{eo}
{\sc S.~Esedo{\=g}lu and F.~Otto}, {\em Threshold dynamics for networks with
  arbitrary surface tensions}, Communications on Pure and Applied Mathematics,
  68 (2015), pp.~808--864, \url{https://doi.org/10.1002/cpa.21527}.

\bibitem{eyre1998unconditionally}
{\sc D.~J. Eyre}, {\em Unconditionally gradient stable time marching the
  {Cahn-Hilliard} equation}, MRS Online Proceedings Library Archive, 529
  (1998).

\bibitem{mbo92}
{\sc B.~Merriman, J.~Bence, and S.~Osher}, {\em Diffusion generated motion by
  mean curvature}, The Computational Crystal Growers. AMS Selection in Math.,
  (1992), pp.~73--83.

\bibitem{mbo94}
{\sc B.~Merriman, J.~K. Bence, and S.~J. Osher}, {\em Motion of multiple
  junctions: a level set approach}, J. Comput. Phys., 112 (1994), pp.~334--363,
  \url{https://doi.org/10.1006/jcph.1994.1105}.

\bibitem{shen2018scalar}
{\sc J.~Shen, J.~Xu, and J.~Yang}, {\em The scalar auxiliary variable (sav)
  approach for gradient flows}, Journal of Computational Physics, 353 (2018),
  pp.~407--416.

\bibitem{shen2010numerical}
{\sc J.~Shen and X.~Yang}, {\em Numerical approximations of {Allen-Cahn} and
  {Cahn-Hilliard} equations}, Discrete Contin. Dyn. Syst, 28 (2010),
  pp.~1669--1691.

\bibitem{shin2017unconditionally}
{\sc J.~Shin, H.~G. Lee, and J.-Y. Lee}, {\em Unconditionally stable methods
  for gradient flow using convex splitting {Runge--Kutta} scheme}, Journal of
  Computational Physics, 347 (2017), pp.~367--381.

\bibitem{shin2020energy}
{\sc J.~Shin and J.-Y. Lee}, {\em An energy stable {Runge--Kutta} method for
  convex gradient problems}, Journal of Computational and Applied Mathematics,
  367 (2020), p.~112455.

\bibitem{alex2019second}
{\sc A.~Zaitzeff, S.~Esedoglu, and K.~Garikipati}, {\em Second order threshold
  dynamics schemes for two phase motion by mean curvature}, 2019,
  \url{https://arxiv.org/abs/1911.05110}.

\end{thebibliography}
\end{document}